\documentclass[journal]{IEEEtran}
\ifCLASSINFOpdf
\else
\fi

\IEEEoverridecommandlockouts
\usepackage{epsfig}
\usepackage{color}
\usepackage{amssymb,latexsym,amsfonts,amsmath,mathrsfs}
\usepackage{graphicx}
\usepackage{cite}
\usepackage{bbm}
\usepackage[linesnumbered,ruled]{algorithm2e}
\allowdisplaybreaks[4]

\newtheorem{thm}{Theorem}
\newtheorem{lem}{Lemma}
\newtheorem{definition}{Definition}
\newtheorem{assum}{Assumption}
\newtheorem{rem}{Remark}

\newcommand{\N}{{\mathbb{N}}}

\begin{document}
%
\title{Distributed Observer for Descriptor Linear System: The Luenberger Observer Method}
%
%
%

\author{Shuai~Liu,
        and~Haotian~Xu
\thanks{This paper is funded by National Natural Science Foundation of China No. 62303273, 62373226, 62133008, Natural Science Foundation of Shandong Province No. ZR2023QF072, and China Postdoctoral Science Foundation No. 2022M721932. The corresponding author: Haotian Xu.}
\thanks{S. Liu and H. Xu are with School of Control Science and Engineering, Shandong University, Jinan 250061~(e-mail: xuhaotian\_1993@126.com; liushuai@sdu.edu.cn)}
}

%
%

\markboth{Submitted to IEEE Transactions on Control of Network Systems}%
{Haotian Xu \MakeLowercase{\textit{et al.}}: Distributed Observer Design}
%



\maketitle

\begin{abstract}
This paper concerns the distributed observer for the descriptor linear system. Unlike centralized descriptor system observers, in the case of distributed observers, each agent either finds it difficult to independently eliminate impulses, or the observer dynamics after eliminating pulses cannot be implemented. To overcome this issue, this paper develops the structure of the distributed observer in two different scenarios, and the observer parameters are presented through a novel design. Moreover, we provide two implementation methods for distributed observer in different scenarios. As a result, each local observer has the ability to reconstruct the states of the underlying system, including its impulse phenomenon. Finally, simulation results verify the validity of our results.
\end{abstract}


%
\IEEEpeerreviewmaketitle

\section{Introduction}\label{sec1}

In recent years, a type of high-performance distributed control law has entered the academic field of vision. \cite{XU2023ISA} and \cite{Xu2022TIV} pointed out that the distributed control law established based on the distributed observer has the ability to approximate the performance of centralized control arbitrarily, and the implementation of this ability mainly benefits from the ability of the state omniscience of distributed observer. Similar studies can also be seen in \cite{Kim2016distributed,Liu2018Cooperative,Xu2020Distributed}, where relevant theories are applied to fields including distributed games \cite{8985536}, distributed power grid control \cite{10543163,2021DistributedMeng}, distributed flexible structure control \cite{Zhang2018Distributed}, attack detection \cite{10543163,SU2024111370}, and so on. Accordingly, distributed observer have received widespread attention and quickly attracted the interest of a large number of scholars \cite{Han2017A,2019Completely,Ortega2022distributed,Silm2020A,9676425,Xu2022SMC}. Different forms of distributed observers, such as distributed function observer \cite{8093658,Martins2017Design}, distributed minimum order observer \cite{Han2018Towards,10400937}, distributed interval observer \cite{Wang2023interval,Huang2022interval,LI2022146}, adaptive distributed observer \cite{8730469}, and distributed observer with unknown input \cite{YANG2022110631,cao2023distributed,cao2023distributedSCL}, have been proposed successively. Some studies focused on distributed observers for different types of systems. For example, \cite{Xu2021IJRNC,battilotti2019distributed,Xu2024JAS} have studied the distributed observer for nonlinear systems independently. \cite{rego2023distributed} and \cite{Bertollo2023hybrid} have developed distributed observers for linear Time-variant systems and linear hybrid systems, respectively. Additionally, distributed observer with event triggered communication \cite{9461598}, network communication delay \cite{8118289,wang2022Dis,mast2023unified}, switching topology network \cite{Xu2022Cybernetics,Xu2025TAC,zhang2023distributed} have also been reported.

All of the above results concern normal system situations. However, many practical physical systems often need to be described by descriptor systems, which have found their
application in modeling the motion of aircraft and in chemical processes,
the mineral industry, electrical circuits, economic systems, and robotics \cite{OSORIOGORDILLO20193,Belov2018,8611289,10243619,10756229}. When facing the distributed control problem of the aforementioned descriptor systems, we also hope to establish a distributed control law with centralized performance. Therefore, it is necessary to first study the distributed observer problem for descriptor systems. Up to now, the observer theory of descriptor systems has developed relatively mature, but research on distributed observer of descriptor systems is still blank.

Though centralized observers for descriptor systems have made mature research progress \cite{rios2024design,tripathi2024observer,Wang2024interval}, it should be pointed out that the distributed observer for descriptor systems is by no means a simple combination of the normal systems' distributed observer theories and the descriptor systems' centralized observer theories. {\color{blue}Specifically, the distributed observer of descriptor systems mainly faces the following challenges. First, the distributed observer of normal systems only needs to consider the stability of the error system \cite{Han2017A,Han2018Towards,2019Completely,Xu2021IJRNC,Xu2024JAS,8093658}, and need not deal with impulse elimination problem, which is the unique requirement of descriptor systems. What is more difficult is that the impulse elimination and the stability analysis of the descriptor systems are interrelated and cannot be achieved independently. Therefore, it brings a challenge to the design of distributed observers while simultaneously considering stability and impulse elimination issues. In addition, in the case of a centralized observer of the descriptor system, the invertible matrix obtained by eliminating impulses can ensure that the observer system is implementable \cite{duan2010analysis}. However, in distributed observers, the implementability of distributed observers can only be guaranteed through traditional methods if the observer parameters are solved in a centralized way. The achievement of descriptor function observer can alleviate the problem of implementation to a certain extent \cite{alma2023adaptive,alma2014adaptive}. However, these methods eliminate impulses based on global information. Therefore, either the observer parameters cannot be solved in a distributed manner, or the local observers cannot be implemented if we design a distributed observer for a descriptor system based on existing centralized observer theories. As a result, another challenge of this paper is to study the implementation methods after designing a distributed observer. }

{\color{cyan}Considering the above issues, this paper contributes to the following three aspects. First, it is the first attempt to propose a distributed observer for the descriptor system. Specifically, we contribute to the design of distributed observers for descriptor systems under two different assumptions, which can effectively estimate all system states, including impulses. In addition, when the descriptor system degenerates into a normal system, the distributed observers in both cases can degenerate into the distributed observers of the normal linear system. Moreover, we have designed implementation methods for distributed observers under two different assumptions, ensuring the normal operation of distributed observers in descriptor systems.}


The remainder of this paper is organized as follows. Section \ref{sec2} introduces some basic knowledge of the descriptor system and formulates the problem. The main results of this paper shown in Section \ref{sec3} and Section \ref{sec4} introduce the implementation methods of the distributed observer. Section \ref{sec5} simulates the developed methods, and Section \ref{sec6} concludes this paper.

\section{Preliminaries and Problem Formulation}\label{sec2}

\textit{Notations:}
(N1)~Let $\mathbb{R}^n$ be the Euclidean space of $n$-dimensional column vectors and $\mathbb{R}^{m\times n}$ be the Euclidean space of $m\times n$-dimensional matrices. $0_{m\times n}\in\mathbb{R}^{m\times n}$ is the zero matrix and $I_n$, $0_n$ stand for identity matrix and all zero square matrix, respectively. In addition, for some zero matrices with dimension that do not produce ambiguity, we will omit their subscripts. Denote $A^\top $ the transpose of matrix $A$. We cast $col\{A_1,\ldots,A_n\}$ as $[A_1^\top ,\ldots,A_n^\top ]^\top $, and $diag\{A_1,\ldots,A_n\}$ as a block diagonal matrix with $A_i$ on its diagonal, where $A_1,\ldots,A_n$ are matrices with arbitrary dimension. $dim\{A\}$ denotes the dimension of square matrix $A$. ${\rm ker} A$ and ${\rm image} A$ represent the kernel space and image space of $A$ respectively. Furthermore, we denote $J_m\in\mathbb{R}^{m\times m}$ the nilpotent matrix. $\|\cdot\|$ is the $2$-norms of vectors or matrices.

(N2)~Set $\mathcal{G}=\{\mathcal{V},\mathcal{E},\mathcal{A}\}$ be a directed graph, $\mathcal{V}$ and $\mathcal{E}$ be the set of its nodes and arcs respectively and denote by $\mathcal{A}$ its adjacency matrix. The element $\alpha_{ij}=1$ of $\mathcal{A}=[\alpha_{ij}]_{i,j=1}^N$ indicates an arc pointing from node $j$ to $i$, and $\alpha_{ij}=0$ otherwise. A directed path from $i_1$ to $i_k$ composed of a set of nodes $\{i_1,\ldots,i_k\}$ is defined by $\alpha_{i_{p+1}i_p}=1$ for all $p\in\{1,\ldots,k\}$. We say the network is strongly connected if there is a path between any pair of nodes $(i,j)$. An input-degree matrix $D=diag\{d_1,\ldots,d_N\}$ is defined by $d_i=\sum_{j=1}^N\alpha_{ij}$, and the matrix $\mathcal{L}=D-\mathcal{A}$ is graph Laplacian of $\mathcal{G}$. The directed graph is a strongly connected graph if there is a directed path between any pair of nodes belonging to $\mathcal{V}$.

\subsection{Preliminaries}
A descriptor linear system often appear in the following form 
\begin{align}
	&E\dot{x}=Ax,\label{sys-11}\\
	&y=Cx,\label{sys-12}
\end{align}
where $x\in\mathbb{R}^n$ and $y\in\mathbb{R}^p$ are the states and outputs of the descriptor system, respectively; $E\in\mathbb{R}^{n\times n}$, $A\in\mathbb{R}^{n\times n}$, and $C\in\mathbb{R}^{p\times n}$ are the system coefficient matrices. In this paper, the system (\ref{sys-11}) is assumed to be regular, i.e., there exists a constant $c$ such that $det\{cE-A\}\neq 0$. In this scene, two non-singular matrices $Q_\star,P_\star\in\mathbb{R}^n$ can be found such that
\begin{align*}
	&Q_\star EP_\star=diag\{I_{n_1},N_2\},\\
	&Q_\star AP_\star=diag\{A_1,I_{n_2}\},\\
	&CP_\star=[C_{\star,1},C_{\star,2}],
\end{align*}
where $A_1\in\mathbb{R}^{n_1\times n_1},~N_2\in\mathbb{R}^{n_2\times n_2}$ with $n_1+n_2=n$ and $N_2$ being a nilpotent matrix. Accordingly, the so-called standard decomposition form (SDF) is introduced as
\begin{align}
	&\dot{x}_1=A_1x_1,\label{sys-21}\\
	&N_2\dot{x}_2=x_2,\label{sys-22}\\
	&y=C_{\star,1}x_1+C_{\star,2}x_2,\label{sys-23}
\end{align}
where $col\{x_1,x_2\}=P_\star^{-1}x$, $C_{\star,1}\in\mathbb{R}^{p\times n_1}$, $C_{\star,2}\in\mathbb{R}^{p\times n_2}$; $x_1\in\mathbb{R}^{n_1},~x_2\in\mathbb{R}^{n_2}$ are the states of slow subsystem and fast subsystem, respectively. 

Furthermore, we can also find two non-singular matrices $Q_\diamond,P_\diamond\in\mathbb{R}^n$ such that
\begin{align*}
	&Q_\diamond EP_\diamond=diag\{I_{l},0_{(n-l)\times(n-l)}\},\\
	&Q_\diamond AP_\diamond=\begin{bmatrix}A_{11}&A_{12}\\A_{21}&A_{22}\end{bmatrix},\\
	&CP_\diamond=[C_{\diamond,1},C_{\diamond,2}],
\end{align*}
where $l$ is the number of non-zero eigenvalues of $E$, and $A_{11}\in\mathbb{R}^{l\times l}$, $A_{12}\in\mathbb{R}^{l\times (n-l)}$, $A_{21}\in\mathbb{R}^{(n-l)\times l}$, $A_{22}\in\mathbb{R}^{(n-l)\times (n-l)}$, $C_{\diamond,1}\in\mathbb{R}^{p\times l}$, $C_{\diamond,2}\in\mathbb{R}^{p\times (n-l)}$. Then, dynamic decomposition form (DDF) with $col\{x_1,x_2\}=P_\diamond^{-1}$ can be expressed as
\begin{align}
	\dot{x}_{1}&=A_{11}x_1+A_{12}x_2,\label{ddf1}\\
	0&=A_{21}x_1+A_{22}x_2,\label{ddf2}\\
	y&=C_{\diamond,1}x_1+C_{\diamond,2}x_2.
\end{align}

\begin{definition}[\cite{duan2010analysis}]\label{def-observable}
	System (\ref{sys-11})--(\ref{sys-12}) is R-observable if and only if the pair $(C_{\star,1},A_1)$ is observable. Furthermore, system (\ref{sys-11}) and (\ref{sys-12}) is I-observable if the impulse can be observed by output measurements. The necessary and sufficient conditions of I-observable is ${\rm ker} N_2\cap {\rm image} N_2\cap {\ker C_{\star,2}}=\{0\}$. We call the system (\ref{sys-11}) and (\ref{sys-12}) is C-observable if the pairs $(C_{\star,1},A_1)$ and $(C_{\star,2},N_2)$ are both observable.
\end{definition}

\begin{definition}[\cite{duan2010analysis}]\label{def-impules}
	There is an impulse term in the dynamic response of the system (\ref{sys-11}) if and only if $N_2\neq 0$, or $A_{22}$ is not invertible. We call the system (\ref{sys-11}) impulse-free if its dynamics do not contain impulse.
\end{definition}

\begin{definition}[\cite{duan2010analysis}]\label{def-admissble}
	System (\ref{sys-11}) is admissible if its state $x$ can converge to zero asymptotically and the dynamic of $x$ is impulse-free for any initial states. For convenience, we also call $(E,A)$ admissible when the system (\ref{sys-11}) is admissible.
\end{definition}

\begin{lem}[\cite{duan2010analysis}]\label{observable}
	System (\ref{sys-11}) and (\ref{sys-12}) is R-observable and I-observable if it is C-observable.
\end{lem}

\begin{lem}[\cite{duan2010analysis}]\label{admissble}
	System (\ref{sys-11}) is admissible if there exists two symmetric matrices $X\geq 0$ and $Y>0$ satisfying 
	\begin{align}\label{lyap}
		E^\top XA+A^\top XE=E^\top YE. 
	\end{align}
	Moreover, there is a symmetric positive definite $X>0$ solved by (\ref{lyap}) for arbitrary $Y>0$ if the system (\ref{sys-11}) is admissible.
\end{lem}

\subsection{Problem formulation}

In this paper, we assume that system (\ref{sys-11}) consists of $N$ agents and the output information is partitioned into $y=col\{y_1,\ldots,y_N\}$ and $C=col\{C_1,\ldots,C_N\}$, where $y_i\in\mathbb{R}^{p_i}$ and $C_i\in\mathbb{R}^{p_i\times n}$ with $\sum_{i=1}^Np_i=p$. The portion $y_i=C_ix$ is assumed to be the only information that can be obtained by agent $i$. $N$ agents can share their information through a communication network $\mathcal{G}$.

Each agent needs to design a local observer and is supposed to reconstruct $x$ of the system (\ref{sys-11}) with the help of its measurement outputs and the information exchanged from their neighbors via a communication network. Let $\hat{x}_i$ be the state estimation generated by the $i$th local observer on agent $i$, and the overall goal of this paper is to achieve
\begin{align}\label{omsci}
	\lim_{t\to\infty}\|\hat{x}_i-x\|=0,
\end{align}
for $i=1,\ldots,N$. All $N$ local observers and the communication network constitute the so-called distributed observer. We say the distributed observer can achieve omniscience asymptotically if (\ref{omsci}) is satisfied. 

To achieve the aforementioned goals, a basic assumptions are necessary, and it is the most common assumption in the field of distributed observer.

\begin{assum}\label{assum1}
	Communication network $\mathcal{G}$ among $N$ agents is assumed to be a strongly connected graph. 
\end{assum}
\begin{lem}[\cite{Han2017A}]\label{lem1}
	Assume $\mathcal{G}$ is a strongly connected directed graph. Then, there exists a unique set of positive constants $r_1,\ldots,r_N$ such that $r^\top \mathcal{L}=0$, $r^\top 1_N=N$ and $\hat{\mathcal{L}}\triangleq R\mathcal{L}+\mathcal{L}^\top R\geq 0$ with $\hat{\mathcal{L}}1_N=1_N\hat{\mathcal{L}}=0$,
	where $r=col\{r_1,\ldots,r_N\}$, $R=diag\{r_1,\ldots,r_N\}$, and $1_N\in\mathbb{R}^N=col\{1,1,\ldots,1\}$.
\end{lem}

\section{Main results}\label{sec3}

Due to the requirement to balance the implementability of the distributed observer and the admissibility of its error dynamics, we cannot design a distributed Luenberger observer for any descriptor system with C-observability. To cover as many systems as possible, this paper will design distributed observers based on two different assumptions. These two different distributed observers will be established based on SDF and DDF in subsection \ref{sec3.1} and \ref{sec3.2}, respectively. Finally, subsection {\ref{sec3.3}} will compare two designs in-depth.

\subsection{Distributed observer with SDF}\label{sec3.1}
In this section, we will show the design method of the distributed observer for the descriptor system (\ref{sys-11}). {\color{cyan}The basic assumption of this subsection is shown as follows, which is the fundamental assumption of observability in the descriptor system.}

\begin{assum}\label{assum2}
	Based on the full output information $y$, system (\ref{sys-11}) and (\ref{sys-12}) is assumed to be R-observable and I-observable. Additionally, based on the local information $y_i=C_ix$ for any $i=1,\ldots,N$, the system is unnecessary to be R-observable.
\end{assum}

In this scenario, the local observer on the $i$th agent takes the following form:
\begin{align}\label{do}
	&E\dot{\hat{x}}_i=A\hat{x}_i+H_i(y_i-C_i\hat{x}_i)+\gamma W_i^{-1}\sum_{j=1}^N\alpha_{ij}\left(E\hat{x}_j-E\hat{x}_i\right),
\end{align}
where $\hat{x}_i\in\mathbb{R}^n$ and $\hat{x}_j\in\mathbb{R}^n$ are the estimated states generated by agent $i$ and $j$ respectively; observer gain matrix $H_i\in\mathbb{R}^{n\times p}$, coupling gain $\gamma$, and weighted matrix $W_i\in\mathbb{R}^{n\times n}$ are parameters that will be designed later. 

{\color{magenta}Before this, it is supposed to focus the SDF to do some preparation work.} We rewrite $C_{\star,1}$, $C_{\star,2}$ as $C_{\star,1}=col\{C_{\star,11},\ldots,C_{\star,N1}\}$, $C_{\star,2}=col\{C_{\star,12},\ldots,C_{\star,N2}\}$ with $C_{\star,i1}\in\mathbb{R}^{p_i\times n_1}$ and $C_{\star,i2}\in\mathbb{R}^{p_i\times n_2}$. If the pairs $(C_{\star,i1},A_1)$ and $(C_{\star,i2},N_2)$ are not observable, we can find orthogonal matrices $T_{\star,i1}$ and $T_{\star,i2}$ such that
\begin{align}
	T_{\star,i1}^\top A_1T_{\star,i1}=&\begin{bmatrix}A_{1,io}&0\\A_{1,ir}&A_{1,iu}\end{bmatrix},\label{A1}\\
	C_{\star,i1}T_{\star,i1}=&[C_{\star,i1o},0_{p_i\times(n-v_{i1})}],\\
	T_{\star,i2}^\top N_2T_{\star,i2}=&\begin{bmatrix}N_{2,io}&0\\N_{2,1r}&N_{2,1u}\end{bmatrix},\\
	C_{\star,i2}T_{\star,i2}=&[C_{\star,i2o},0_{p_i\times(n-v_{i2})}],\label{A2}
\end{align}
where $A_{1,io}\in\mathbb{R}^{v_{i1}\times v_{i1}}$, $N_{2,io}\in\mathbb{R}^{v_{i2}\times v_{i2}}$; $v_{i1}$ and $v_{i2}$ are the observability index of $(C_{\star,i1},A_1)$ and $(C_{\star, i2},N_2)$, respectively; $A_{1,ir}$, $A_{1,iu}$, $N_{2,ir}$, and $N_{2,iu}$ are matrices with compatible dimension. As a result, according to Lemma \ref{observable}, the observable subsystem 
\begin{align*}
	&\dot{x}_{1,io}=A_{1,io}x_{1,io},\\
	&N_{2,io}\dot{x}_{2,io}=x_{2,io},
\end{align*}
is C-observable. 

According to \cite{Han2017A}, we can obtain the following lemma by constructing $T_{\star i}=diag\{T_{\star,i1},T_{\star,i2}\}$ (the proof is the same as \cite[Lemma 2]{Han2017A}).
\begin{lem}\label{jointob}
	Let $T_\star=diag\{T_{\star 1},\ldots,T_{\star N}\}$ and construct $G_\star=diag\{G_{\star 1},\ldots,G_{\star N}\}$, where $G_{i\star}=diag\{G_{\star,i1},G_{\star,i2}\}$ for $i=1,\ldots,N$ with $G_{\star,i1}=diag\{g_iI_{v_{i1}},0_{(n-v_{i1})\times (n-v_{i1})}\}$ and $G_{\star,i2}=diag\{g_iI_{v_{i2}},0_{(n-v_{i2})\times (n-v_{i2})}\}$. Then, for arbitrary $g_i>0$, there exists $\mu>0$ such that
	\begin{align}
		T_\star^\top \left(\mathcal{L}\otimes I_n\right)T_\star+G_\star>\mu I_{nN}.
	\end{align}
\end{lem}

{\color{magenta}Now, based on the aforementioned observable decomposition and Lemma, the parameter design methods and the error dynamics stability of (\ref{do}) will be described by the following theorem.}

\begin{thm}\label{thm1}
	Consider the descriptor system (\ref{sys-11})--(\ref{sys-12}) subject to Assumption \ref{assum2} and assume that the communication network satisfies Assumption \ref{assum1}. Then, the sufficient conditions of distributed observer (\ref{do}) to achieve omniscience asymptotically includes:\\
	{\color{cyan}1)~Observability indexes satisfy $v_{i2}=n_2$ for all $i=1,\ldots,N$ and coupling gain $\gamma>\bar{\sigma}(\tilde{Y}_u)/\mu$, where $\tilde{Y}_{u}=diag\{\tilde{Y}_{1,u},\ldots,\tilde{Y}_{N,u}\}$ and
	\begin{align*}
		\tilde{Y}_{i,u}=diag\left\{r_i\begin{bmatrix}
			&A_{1,ir}^\top \\A_{1,ir}&sym\{A_{1,iu}\}
		\end{bmatrix},~i=1,\ldots,N\right\};
	\end{align*}}\\
	2)~Observer gain is chosen as $H_i=Q_\star^{-1}col\left\{H_{i1},~H_{i2}\right\}$, where 
	\begin{align}
		H_{i1}=T_{\star,i1}\begin{bmatrix}H_{i1o}\\0_{(n_1-v_{i1})\times p_i}\end{bmatrix},~H_{i2}=T_{\star,i2}\begin{bmatrix}H_{i2o}\\0_{(n_2-v_{i2})\times p_i}\end{bmatrix},
	\end{align}
	and $H_{i1o}$, $H_{i2o}$ are designed such that $(\tilde{E}_{i,o},\tilde{A}_{i,o}-H_{i,o}C_{\star,io})$ is admissible, where
	\begin{align}
		&\tilde{E}_{i,o}=\begin{bmatrix}
			I_{v_{i1}}&\\&N_{2,io}
		\end{bmatrix},~\tilde{A}_{i,o}=\begin{bmatrix}
			A_{1,io}&\\&I_{v_{i2}}
		\end{bmatrix},\\
		&H_{i,o}=col\{H_{i1o},H_{i2o}\},~C_{\star,io}=[C_{\star,i1o},C_{\star,i2o}];
	\end{align}
	3)~Weighted matrix $W_i$ takes the form of:
	\begin{align}
		&W_i=Q_\star^{-1}W_{Ti}Q_\star=Q_\star^{-1}\begin{bmatrix}W_{T,i11}&W_{T,i12}\\W_{T,i21}&W_{T,i22}\end{bmatrix}Q_\star,\\
		&W_{Ti}=T_{\star i}\tilde{W}_iT_{\star i}^\top ,~i=1,\dots,N,\\
		&\tilde{W}_i=\begin{bmatrix}
			\tilde{W}_{i,11o}&0&\tilde{W}_{i,12o}&0\\
			0&I_{n_1-v_{i1}}&0&0\\
			\tilde{W}_{i,21o}&0&\tilde{W}_{i,22o}&0\\
			0&0&0&I_{n_2-v_{i2}}
		\end{bmatrix},
	\end{align}
	and
	\begin{align}
		\tilde{W}_{i,o}=\begin{bmatrix}
			\tilde{W}_{i,11o}&\tilde{W}_{i,12o}\\
			\tilde{W}_{i,21o}&\tilde{W}_{i,22o}
		\end{bmatrix}
	\end{align}
	is the solution of Lyapunov function
	\begin{align}\label{LYAPF}
		sym\left\{\tilde{E}_{i,o}^\top \tilde{W}_{i,o}\left(\tilde{A}_{i,o}-H_{i,o}C_{\star,io}\right)\right\}=\tilde{E}_{i,o}^\top Y\tilde{E}_{i,o}
	\end{align}
	for arbitrary $Y>0$.
\end{thm}

\begin{IEEEproof}
Let $e_i=\hat{x}_i-x$, and we have
\begin{align}\label{error1}
	E\dot{e}_i=Ae_i-H_iC_ie_i-\gamma W_i^{-1}\sum_{j=1}^Nl_{ij}Ee_j,
\end{align}
where $l_{ij}$ is the element located at row $i$ and column $j$ of Laplacian matrix  $\mathcal{L}$. Based on the theory of SDF, there exists a couple of non-singular matrices $P_\star$ and $Q_\star$ such that $Q_\star EP_\star=diag\{I_{n_1},N_{2}\}$, $Q_\star AP_\star=col\{A_1,I_{n_2}\}$, $C_iP_\star=[C_{\star,i1},C_{\star,i2}]$, where $A_1\in\mathbb{R}^{n_1\times n_1}$, $C_{\star,i1}\in\mathbb{R}^{p_i\times n_1}$, and $C_{\star,i2}\in\mathbb{R}^{p_i\times n_2}$. Accordingly, by setting $\eta_i=col\{\eta_{i1},\eta_{i2}\}=P_\star^{-1}e_i$ and $col\{H_{i1},H_{i2}\}=Q_\star H_i$ with $\eta_{i1}\in\mathbb{R}^{n_1}$, $\eta_{i2}\in\mathbb{R}^{n_2}$, $H_{i1}\in\mathbb{R}^{n_1\times p_i}$, and $H_{i2}\in\mathbb{R}^{n_2\times p_i}$, error dynamics (\ref{error1}) can be transformed into
\begin{align}
	\dot{\eta}_{i1}=&A_1\eta_{i1}-H_{i1}C_{\star,i1}\eta_{i1}-H_{i1}C_{\star,i2}\eta_{i2}\notag\\
	&-\gamma W_{T,i11}^{inv}\sum_{j=1}^Nl_{ij}\eta_{i1}-\gamma W_{T,i12}^{inv}\sum_{j=1}^Nl_{ij}N_2\eta_{i2},\label{error-21}\\
	N_{2}\dot{\eta}_{i2}=&\eta_{i2}-H_{i2}C_{\star,i1}\eta_{i1}-H_{i2}C_{\star,i2}\eta_{i2}\notag\\
	&-\gamma W_{T,i21}^{inv}\sum_{j=1}^Nl_{ij}\eta_{i1}-\gamma W_{T,i22}^{inv}\sum_{j=1}^Nl_{ij}N_{2}\eta_{i2},\label{error-21}
\end{align}
where
\begin{align*}
	\begin{bmatrix}
		W_{T,i11}^{inv}&W_{T,i12}^{inv}\\W_{T,i21}^{inv}&W_{T,i22}^{inv}
	\end{bmatrix}\triangleq \begin{bmatrix}W_{T,i11}&W_{T,i12}\\W_{T,i21}&W_{T,i22}\end{bmatrix}^{-1}.
\end{align*}
By denoting $\eta=col\{\eta_{1},\eta_2,\ldots,\eta_N\}$ and $\eta_i=col\{\eta_{i1},\eta_{i2}\}$ for all $i=1,\ldots,N$, the following compact form is obtained
\begin{align}\label{error-11}
	(I_N\otimes \bar{E})\dot{\eta}=\bar{A}\eta-\gamma \bar{W}^{-1}(\mathcal{L}\otimes I_n)(I_N\otimes \bar{E})\eta,
\end{align}
where 
\begin{align*}
	&\bar{E}=diag\{I_{n_1},N_{2}\},\\
	&\bar{A}=diag\{\bar{A}_1,\ldots,\bar{A}_N\},\\
	&\bar{A}_i=\begin{bmatrix}A_1&\\&I_{n_2}\end{bmatrix}-\begin{bmatrix}H_{i1}\\H_{i2}\end{bmatrix}\begin{bmatrix}C_{\star,i1}&C_{\star,i2}\end{bmatrix},\\
	&\bar{W}=diag\{W_{T1},\ldots,W_{TN}\}.
\end{align*}

{\color{red}Now, by letting $\xi=T_\star^\top \eta$, and multiplying $T_\star$ on both side of (\ref{error-11}), we obtain}
\begin{align}
	T_\star^\top (I_N\otimes \bar{E})T_\star\dot{\xi}=&diag\{T_{1\star}^\top \bar{A}_1T_{1\star},\ldots,T_{N\star}^\top \bar{A}_NT_{N\star}\}\xi\notag\\
	&-\gamma\tilde{W}^{-1}T_\star^\top (\mathcal{L}\otimes I_n)T_\star T_\star^\top (I_N\otimes \bar{E})T_\star\xi.\label{error2}
\end{align}
{\color{red}According to the observability decomposition of SDF mentioned at the beginning of this subsection, we have }
\begin{align*}
	T_\star^\top (I_N\otimes \bar{E})T_\star=&diag\{T_{1\star}^\top ET_{1\star},\ldots,T_{N\star}^\top ET_{N\star}\}\notag\\
	=&diag\left\{\tilde{E}_1,\ldots,\tilde{E}_N\right\},
\end{align*}
where
\begin{align}
	\tilde{E}_i=\begin{bmatrix}I_{v_{i1}}&&&\\&I_{n_1-v_{i1}}&&\\&&N_{2,io}&0\\&&N_{2,ir}&N_{2,iu}\end{bmatrix},
\end{align}
and
\begin{align}
	&T_{\star i}^\top \bar{A}_{i}T_{\star i}\triangleq\tilde{A}_i\notag\\
	=&\begin{bmatrix}\Psi_{1,i1o}&0&-H_{i1o}C_{\star,i2o}&0\\A_{1,i1r}&A_{\star,i1u}&0&0\\-H_{i2o}C_{\star,i1o}&0&\Psi_{v_{i2}}&0\\0&0&0&I_{n_2-v_{i2}}\end{bmatrix},
\end{align}
where $\Psi_{1,i1o}=A_{1,i1o}-H_{i1o}C_{\star,i1o}$, $\Psi_{v_{i2}}=I_{v_{i2}}-H_{i2o}C_{\star,i2o}$.
Hence, we know from the observability decomposition theory of descriptor linear system that the pair $(\tilde{E}_{i,o},\tilde{A}_{i,o},C_{\star,io})$ is C-observable. Therefore, there exists proper $H_{i,o}$ such that $(\tilde{E}_{i,o},\tilde{A}_{i,o}-H_{i,o}C_{\star,io})$ is admissible. Therefore, according to Lemma \ref{admissble}, by giving $\tilde{Y}_{i,o}=\gamma diag\{I_{v_{i1}},I_{v_{i2}}\}$, there exists $\tilde{W}_{i,o}>0$ solved by
\begin{align}\label{lyp11}
	&sym\left\{\tilde{E}_{i,o}^\top \tilde{W}_{i,o}\left(\tilde{A}_{i,o}-H_{i,o}C_{\star,io}\right)\right\}\notag\\
	=&-\tilde{E}_{i,o}^\top \tilde{Y}_{i,o}\tilde{E}_{i,o}=diag\left\{-\gamma I_{v_{i1}},-\gamma N_{2,io}^\top N_{i,2o}\right\}.
\end{align}

{\color{red}To move on, we define several compact forms $\tilde{E}=diag\{\tilde{E}_1,\ldots,\tilde{E}_N\}$, $\tilde{W}=diag\{\tilde{W}_1,\ldots,\tilde{W}_N\}$, and $\tilde{A}=diag\{\tilde{A}_1,\ldots,\tilde{A}_N\}$. The Lyapunov method of descriptor system introduced in Lemma \ref{admissble} can be considered for (\ref{error2}) by using $\tilde{E}$, $\tilde{A}$, and $(R\otimes I_n)\tilde{W}$ and deduce
\begin{align}\label{Lyap1}
	&sym\left\{\tilde{E}^\top (R\otimes I_n)\tilde{W}\left(\tilde{A}-\gamma\tilde{W}^{-1}T_\star^\top (\mathcal{L}\otimes I_n)T_\star\tilde{E}\right)\right\}\notag\\
	=&(R\otimes I_n)\left(\tilde{E}^\top \tilde{W}\tilde{A}+\tilde{A}^\top \tilde{W}\tilde{E}\right)\notag\\
	&-\gamma\tilde{E}^\top T_\star^\top \left(\left(R\mathcal{L}+\mathcal{L}^\top R\right)\otimes I_n\right)T_\star\tilde{E}\notag\\
	=&diag\left\{r_i\begin{bmatrix}
		-\gamma I_{v_{i1}}&A_{1,ir}^\top \\
		A_{1,ir}&sym\{A_{1,iu}\}\end{bmatrix},\right.\notag\\
	 &\quad\quad~\left.r_i\begin{bmatrix}
		-\gamma N_{2,io}^\top N_{i,2o}&N_{2,ir}^\top \\
		N_{2,ir}&sym\{N_{2,iu}\}
	\end{bmatrix},~i=1,\ldots,N\right\}\notag\\
	&-\gamma\tilde{E}^\top T_\star^\top \left(\left(R\mathcal{L}+\mathcal{L}^\top R\right)\otimes I_n\right)T_\star\tilde{E}.
\end{align} 
}Note that---for all $i=1,\ldots,N$---if and only if $v_{i2}=n_2$, there are symmetric positive definite matrix $\tilde{Y}_i=-\gamma diag\{r_iI_{v_{i1}},0_{(n-v_{i1})\times(n-v_{i1})},r_iI_{n_2}\}$ and {\color{red}symmetric matrix $\tilde{Y}_{iu}$ such that
\begin{align}
	&\tilde{E}^\top \tilde{W}\tilde{A}+\tilde{A}^\top \tilde{W}\tilde{E}\notag\\
	=&\begin{bmatrix}
		-\gamma I_{v_{i1}}&A_{1,ir}^\top &&\\
		A_{1,ir}&sym\{A_{1,iu}\}&&\\
		&&-\gamma N_{2,io}^\top N_{i,2o}&N_{2,ir}^\top \\
		&&N_{2,ir}&sym\{N_{2,iu}\}
	\end{bmatrix}\notag\\
	=&\begin{bmatrix}
	-\gamma I_{v_{i1}}&A_{1,ir}^\top &&\\
	A_{1,ir}&sym\{A_{1,iu}\}&&\\
	&&-\gamma N_{2,io}^\top N_{i,2o}
\end{bmatrix}\notag\\
	=&\tilde{E}_{i}^\top \tilde{Y}_i\tilde{E}_i+\tilde{E}_i^\top \tilde{Y}_{i,u}\tilde{E}_i.\label{Y}
\end{align}
Herein, the second equal sign is because $v_{i2}=n_2$ which leads thereby to $dim\{N_{2,iu}\}=0$.} Now, we construct $\tilde{Y}_o=diag\{\tilde{Y}_{1},\ldots,\tilde{Y}_{N}\}$. Then, (\ref{Lyap1}) can be rewritten as
\begin{align}
	&sym\left\{\tilde{E}^\top (R\otimes I_n)\tilde{W}\left(\tilde{A}-\gamma\tilde{W}^{-1}T_\star^\top (\mathcal{L}\otimes I_n)T_\star\tilde{E}\right)\right\}\notag\\
	=&\tilde{E}^\top \tilde{Y}_{o}\tilde{E}+\tilde{E}^\top \tilde{Y}_{u}\tilde{E}-\gamma\tilde{E}^\top T_\star^\top \left(\left(R\mathcal{L}+\mathcal{L}^\top R\right)\otimes I_n\right)T_\star\tilde{E}\notag\\
	=&\tilde{E}^\top \left(\tilde{Y}_{o}-\gamma T_\star^\top \left(\hat{\mathcal{L}}\otimes I_n\right)T_\star+\tilde{Y}_u\right)\tilde{E}.
\end{align}

{\color{red}In light of Lemma \ref{jointob}, we know there exists a constant $\mu>0$ such that
\begin{align}
	\tilde{Y}_{o}-\gamma T_\star^\top \left(\hat{\mathcal{L}}\otimes I_n\right)T_\star<-\gamma\mu I_{nN}.
\end{align}
Therefore, according to Lemma \ref{admissble}, error dynamics (\ref{error-11}) is admissible if there is a $\gamma$ such that  
\begin{align}
	&\tilde{Y}_{o}-\gamma T_\star^\top \left(\hat{\mathcal{L}}\otimes I_n\right)T_\star+\tilde{Y}_u=-\gamma\mu I_{nN}+\tilde{Y}_u<0.
\end{align}
The sufficient condition of this inequality is $\gamma>\bar{\sigma}(\tilde{Y}_u)/\mu$.}

Consequently, we have proved that $\lim_{t\to\infty}e_i(t)=0$ and the dynamics of $e_i$ are impulse-free for all $i=1,\ldots,N$ if conditions described in Theorem \ref{thm1} are fulfilled. 
\end{IEEEproof}

This section gives the design method of distributed observer for descriptor linear system with the assumption $v_{i2}=n_2$ and Assumption \ref{assum2}, and describes the design method of parameters in Theorem \ref{thm1} in detail. {\color{blue}For the convenience of readers using this method, we provide a brief guideline in follows for designing.\\
\textbf{STEP 1:}~Let $\hat{E}=(cE-A)^{-1}E$ and calculate 
\begin{align*}
	Q_\star=diag\{\hat{E}_1^{-1},(c\hat{E}_2-I)^{-1}\}U(cE-A)^{-1},~P_\star=U^{-1},
\end{align*}
where $U$ is a non-singular matrix such that $\hat{E}$ can be transformed into Jordan canonical form $U\hat{E}U^{-1}=diag\{\hat{E}_1,\hat{E}_2\}$ with $\hat{E}_1$ being a non-singular matrix and $\hat{E}_2$ being a nilpotent matrix. Then, we can use $Q_\star$ and $P_{star}$ to obtain SDF for descriptor system.\\
\textbf{STEP 2:}~Calculate matrices $A_{1,io}$, $A_{1,ir}$, and so on based on the observability decomposition (\ref{A1})--(\ref{A2}) of SDF.\\
\textbf{STEP 3:}~Calculate observer gain $H_i$ based on pole placement such that $(\tilde{E}_{i,o},\tilde{A}_{i,o}-H_{i,o}C_{\star,io})$ is admissible. (The pole placement method for descriptor systems can be found in many textbooks on descriptor linear systems, such as \cite{duan2010analysis}.)\\
\textbf{STEP 4:}~For any given $Y>0$, solve weight matrix $\tilde{W}_{i,o}$ from Lyapunov function (\ref{LYAPF}) (Solution for descriptor Lyapunov function can be obtained directly by MATLAB).\\
\textbf{STEP 5:}~Coupling gain $\gamma$ can be chosen by $\gamma>\bar{\sigma}(\tilde{Y}_u)/\mu$, or chosen by adaptive law.}

\begin{rem}
	\color{cyan}Note that calculating $\gamma$ in Theorem \ref{thm1} requires the eigenvalues of Laplacian matrix $\hat{\mathcal{L}}$, which is a global information. To overcome this problem such that the developed method is completely distributed, some methods for distributed estimate the eigenvalues of $\hat{\mathcal{L}}$ can be introduced \cite{TRAN201556}. Alternatively, a more convenient method is to adopt adaptive coupling gain $\gamma=\omega$ and $\dot{\omega}=\|\sum_{j=1}^N\alpha_{ij}(E\hat{x}_j-E\hat{x}_i)\|^2$. 
    The proof of adaptive coupling gain is already mature in distributed observers of normal systems \cite{2019Completely,Xu2025TAC}, and the related derivation process is also applicable in descriptor systems, so we do not repeat it here.
\end{rem}

\subsection{Distributed observer with DDF}\label{sec3.2}

Distributed observer cannot be designed with SDF if condition $v_{i2}=n_2$ required in subsection \ref{sec3.1} is not satisfied. To make up for this deficiency, this subsection studies distributed observer with DDF in order to not rely on $v_{i2}=n_2$.

To this end, we first transform the descriptor system (\ref{sys-11})--(\ref{sys-12}) into DDF, which is shown in (\ref{ddf1})--(\ref{ddf2}). Then, for saving of convenience, we rewrite $C_{\diamond,1}$, $C_{\diamond,2}$ {\color{magenta}(the symbols obtained by DDF)} as $C_{\diamond,1}=col\{C_{\diamond,11},\ldots,C_{\diamond,N1}\}$, $C_{\diamond,2}=col\{C_{\diamond,12},\ldots,C_{\diamond,N2}\}$ with $C_{\diamond,i1}\in\mathbb{R}^{p_i\times l}$ and $C_{\diamond,i2}\in\mathbb{R}^{p_i\times (n-l)}$. The basic assumption of this subsection is given below (A detailed interpretation of this assumption can be found in Remark \ref{rem-assum3}).
\begin{assum}\label{assum3}
	Based on the full output information $y$, the pairs $(C_{\diamond 1},A_{11})$ and $(C_{\diamond 2}, A_{22})$ are observable. In addition, based on the local information $y_i=C_{\diamond,i}x$ for any $i=1,\ldots,N$, the pairs $(C_{\diamond, i1},A_{11})$ and $(C_{\diamond, i2}, A_{22})$ are unnecessary to be observable. 
\end{assum}

In this scenario, the local observer for the $i$th agent can be designed as
\begin{align}\label{do2}
	E\dot{\hat{x}}_i=A\hat{x}_i+H_i(y_i-C_i\hat{x}_i)-\gamma W_i^{-1}\sum_{j=1}^N\alpha_{ij}\left(\hat{x}_j-\hat{x}_i\right),
\end{align}
where $\hat{x}_i\in\mathbb{R}^n$ and $\hat{x}_j\in\mathbb{R}^n$ are with the same definitions as that in (\ref{do}); observer gain matrix $H_i\in\mathbb{R}^{n\times p}$, coupling gain $\gamma$, and weighted matrix $W_i\in\mathbb{R}^{n\times n}$ are parameters that will be designed in Theorem \ref{thm2}.

Before showing the main theorem of this subsection, some notations should be introduced. According to Assumption \ref{assum3}, we can find orthogonal matrices $T_{\diamond,i1}$ and $T_{\diamond,i2}$ such that
\begin{align}
	T_{\diamond,i1}^\top A_{11}T_{\diamond,i1}=&\begin{bmatrix}
		A_{11,io}&0\\A_{11,ir}&A_{11,iu}
	\end{bmatrix},\label{A3}\\
	C_{\diamond,i1}T_{\diamond,i1}=&[C_{\diamond,i1o},0_{p_i\times(l-v_{i1})}],\\
	T_{\diamond,i2}^\top A_{22}T_{\diamond,i2}=&\begin{bmatrix}
		A_{22,io}&0\\A_{22,ir}&A_{22,iu}
	\end{bmatrix},\\
	C_{\diamond,i2}T_{\diamond,i2}=&[C_{\diamond,i2o},0_{p_i\times(n-l-v_{i2})}],\label{A4}
\end{align}
for any $i=1,\ldots,N$, where $v_{i1}$ and $v_{i2}$ are observability indexes of the pair $(C_{\diamond, i1},A_{11})$ and $(C_{\diamond, i2}, A_{22})$, respectively.

\begin{lem}\label{jointob2}
	Let $T_{\diamond,j}=diag\{T_{\diamond, 1j},\ldots,T_{\diamond, Nj}\}$, and construct $G_{\diamond,j}=diag\{G_{\diamond, 1j},\ldots,G_{\diamond, Nj}\}$, for $j=1,2$, where $G_{\diamond,i1}=diag\{g_iI_{v_{i1}},0_{(n-v_{i1})\times (n-v_{i1})}\}$ and $G_{\diamond,i2}=diag\{g_iI_{v_{i2}},0_{(n-v_{i2})\times (n-v_{i2})}\}$. Then, for arbitrary $g_i>0$, there exists $\mu>0$ such that
	\begin{align}
		&T_{\diamond,1}^\top \left(\mathcal{L}\otimes I_l\right)T_{\diamond,1}+G_{\diamond,1}>\mu I_{nN},\\
		&T_{\diamond,2}^\top \left(\mathcal{L}\otimes I_{n-l}\right)T_{\diamond,2}+G_{\diamond,2}>\mu I_{nN}.
	\end{align}
\end{lem}

\begin{thm}\label{thm2}
	Consider the descriptor system (\ref{sys-11})--(\ref{sys-12}) subject to Assumption \ref{assum3}, and the communication network satisfies Assumption \ref{assum1}. Then, the sufficient conditions for distributed observer (\ref{do2}) to achieve omniscience asymptotically include:\\
	1)~Gain matrices are
	\begin{align*}
		&H_i=Q_{\diamond}^{-1}col\{H_{i1},H_{i2}\},\\
		&H_{i1}=T_{\diamond,i1}col\{H_{i1o},0_{(l-v_{i1})\times 1}\},\\
		&H_{i2}=T_{\diamond,i2}col\{H_{i2o},0_{(n-l-v_{i2})\times 1}\},
	\end{align*}
	where $H_{i1o}$ and $H_{i2o}$ are designed such that $A_{11,io}-H_{i1o}C_{\diamond,i1o}$ and $A_{22,io}-H_{i2o}C_{\diamond,i2o}$ are Hurwitz matrices;\\
	2)~Weighted matrix $W_i=P_{\diamond}^{-1}diag\{W_{T,i1},W_{T,i2}\}Q_{\diamond}^{-1}$ and $W_{T,i1}=T_{\diamond,i1}\tilde{W}_{T,i1}T_{\diamond,i1}^\top $, $W_{T,i2}=T_{\diamond,i2}\tilde{W}_{T,i2}T_{\diamond,i2}^\top $, where $\tilde{W}_{T,i1}$ and $\tilde{W}_{T,i2}$ are symmetric positive definite solutions to the following equations
	\begin{align}
		&sym\{\tilde{W}_{T,i1}(A_{11,io}-H_{i1o}C_{\diamond,i1o})\}=-2\gamma I_{v_{i1}},\label{w1}\\
		&sym\{\tilde{W}_{T,i2}(A_{22,io}-H_{i2o}C_{\diamond,i2o})\}=-2\gamma I_{v_{i2}},\label{w2}
	\end{align}
	in which $\mu>0$ is an arbitrary given constant;\\
	3)~Coupling gain $\gamma>0$ satisfies $\gamma>\frac{1}{\mu}\max\{\lambda_1,\lambda_2\}$, where $\lambda_1$ and $\lambda_2$ are given behind (\ref{lambda1}) and (\ref{lyp25}), respectively.
\end{thm}
\begin{IEEEproof} Let $e_i=\hat{x}_i-x$, then we have
\begin{align}\label{e1}
	E\dot{e}_i=(A-H_iC_i)e_i-\gamma W_i^{-1}\sum_{j=1}^N\alpha_{ij}(e_j-e_i).
\end{align}
Bearing in mind the definition of DDF (\ref{ddf1})--(\ref{ddf2}), we know there exists $P_\diamond, Q_\diamond\in\mathbb{R}^n$ such that (\ref{e1}) can be transformed into
\begin{align}
	\dot{\eta}_{i1}=&(A_{11}-H_{i1}C_{\diamond, i1})\eta_{i1}+(A_{12}-H_{i1}C_{\diamond, i2})\eta_{i2}\notag\\
	-&\gamma W_{T,i1}^{-1}\sum_{j=1}^N \alpha_{ij}(\eta_{j1}-\eta_{i1}),\\
	0=&(A_{21}-H_{i2}C_{\diamond, i1})\eta_{i1}+(A_{22}-H_{i2}C_{\diamond, i2})\eta_{i2}\notag\\
	-&\gamma W_{T,i2}^{-1}\sum_{j=1}^N \alpha_{ij}(\eta_{j2}-\eta_{i2}),
\end{align}
where $col\{\eta_{i1},\eta_{i2}\}\triangleq\eta_i$ and $\eta_i=P_\diamond^{-1}e_i$; $col\{H_{i1},H_{i2}\}=Q_\diamond H_i$ with $\eta_{i1}\in\mathbb{R}^{l}$, $\eta_{i2}\in\mathbb{R}^{n-l}$, $H_{i1}\in\mathbb{R}^{l\times p_i}$, and $H_{i2}\in\mathbb{R}^{(n-l)\times p_i}$. 

To move on, we construct two compact forms $\eta_{\diamond 1}=col\{\eta_{11},\ldots,\eta_{N1}\}$ and $\eta_{\diamond 2}=col\{\eta_{12},\ldots,\eta_{N2}\}$, and denote $\Lambda_{ab}=diag\{A_{ab}-H_{ia}C_{\diamond, ib},~i=1,\ldots,N\}$ for $a,b=1,2$, which yields
\begin{align}
	\dot{\eta}_{\diamond 1}=&\Lambda_{11}\eta_{\diamond 1}+\Lambda_{12}\eta_{\diamond 2}-\gamma W_{T,1}^{-1}(\mathcal{L}\otimes I_l)\eta_{\diamond 1},\\
	0=&\Lambda_{21}\eta_{\diamond 1}+\Lambda_{22}\eta_{\diamond 2}-\gamma W_{T,2}^{-1}(\mathcal{L}\otimes I_{n-l})\eta_{\diamond 2},
\end{align}
where $W_{T,1}=diag\{W_{T,i1},~i=1,\ldots,N\}$ and $W_{T,2}=diag\{W_{T,i2},~i=1,\ldots,N\}$. Let $T_{\diamond, 1}=diag\{T_{\diamond,i1},~i=1,\ldots,N\}$, $T_{\diamond, 2}=diag\{T_{\diamond,i2},~i=1,\ldots,N\}$, and $\xi_{\diamond 1}=T_{\diamond 1}^\top \eta_{\diamond 1}$, $\xi_{\diamond 2}=T_{\diamond 2}^\top \eta_{\diamond 2}$. Then, based on conditions 2) and 3), we have
\begin{align}
	\dot{\xi}_{\diamond 1}=&T_{\diamond, 1}^\top \Lambda_{11}T_{\diamond, 1}\xi_{\diamond 1}+T_{\diamond, 1}^\top \Lambda_{12}T_{\diamond, 2}\xi_{\diamond 2}\notag\\
	&-\gamma \tilde{W}_{T,1}^{-1}T_{\diamond, 1}^\top (\mathcal{L}\otimes I_l)T_{\diamond, 1}\xi_{\diamond 1}\notag\\
	\triangleq&\tilde{\Lambda}_{11}\xi_{\diamond 1}+\tilde{\Lambda}_{12}\xi_{\diamond 2},\label{xi-21}\\
	0=&T_{\diamond, 2}^\top \Lambda_{21}T_{\diamond, 1}\xi_{\diamond 1}+T_{\diamond, 2}^\top \Lambda_{22}T_{\diamond, 2}\xi_{\diamond 2}\notag\\
	&-\gamma\tilde{W}_{T,2}^{-1} T_{\diamond, 2}^\top (\mathcal{L}\otimes I_{n-l})T_{\diamond, 2}\xi_{\diamond 2}\notag\\
	\triangleq&\tilde{\Lambda}_{21}\xi_{\diamond 1}+\tilde{\Lambda}_{22}\xi_{\diamond 2},\label{xi-22}
\end{align}
where 
\begin{align*}
	&\tilde{\Lambda}_{11}=T_{\diamond, 1}^\top \Lambda_{11}T_{\diamond, 1}-\gamma \tilde{W}_{T,1}^{-1}T_{\diamond, 1}^\top (\mathcal{L}\otimes I_l)T_{\diamond, 1},\\
	&\tilde{\Lambda}_{12}=T_{\diamond, 1}^\top \Lambda_{12}T_{\diamond, 2},\\
	&\tilde{\Lambda}_{21}=T_{\diamond, 2}^\top \Lambda_{21}T_{\diamond, 1},\\
	&\tilde{\Lambda}_{22}=T_{\diamond, 2}^\top \Lambda_{22}T_{\diamond, 2}-\gamma\tilde{W}_{T,2}^{-1} T_{\diamond, 2}^\top (\mathcal{L}\otimes I_{n-l})T_{\diamond, 2}.
\end{align*}

{\color{red}In what follows, we will prove the admissibility of error dynamics (\ref{xi-21})--(\ref{xi-22}). In light of Lemma \ref{admissble}, we know (\ref{xi-21})--(\ref{xi-22}) is admissible if there are }
\begin{align*}
	X=\begin{bmatrix}X_{11}&X_{12}\\X_{12}^\top &X_{22}\end{bmatrix}\geq 0,~Y=\begin{bmatrix}Y_{11}&Y_{12}\\Y_{12}^\top &Y_{22}\end{bmatrix}> 0.
\end{align*}
such that
\begin{align}
	&\begin{bmatrix}I_{l}&\\&0_{(n-l)\times(n-l)}\end{bmatrix}\begin{bmatrix}X_{11}&X_{12}\\X_{12}^\top &X_{22}\end{bmatrix}\begin{bmatrix}\tilde{\Lambda}_{11}&\tilde{\Lambda}_{12}\\\tilde{\Lambda}_{21}&\tilde{\Lambda}_{22}\end{bmatrix}\notag\\
	+&\begin{bmatrix}\tilde{\Lambda}_{11}^\top &\tilde{\Lambda}_{12}^\top \\\tilde{\Lambda}_{21}^\top &\tilde{\Lambda}_{22}^\top \end{bmatrix}\begin{bmatrix}X_{11}&X_{12}\\X_{12}^\top &X_{22}\end{bmatrix}\begin{bmatrix}I_{l}&\\&0_{(n-l)\times(n-l)}\end{bmatrix}\notag\\
	=&\begin{bmatrix}
		sym\{X_{11}\tilde{\Lambda}_{11}\}+sym\{X_{12}\tilde{\Lambda}_{21}\}&X_{11}\tilde{\Lambda}_{12}+X_{12}\tilde{\Lambda}_{22}\\\tilde{\Lambda}_{12}^\top X_{11}+\tilde{\Lambda}_{22}^\top X_{12}&0
	\end{bmatrix}\notag\\
	=&\begin{bmatrix}I_{l}&\\&0_{(n-l)\times(n-l)}\end{bmatrix}Y\begin{bmatrix}I_{l}&\\&0_{(n-l)\times(n-l)}\end{bmatrix}\notag\\
	=&-diag\{Y_{11},0\}.\label{lyap2}
\end{align}
The above equation holds if and only if
\begin{align}
	&X_{11}\tilde{\Lambda}_{11}+\tilde{\Lambda}_{11}^\top X_{11}+X_{12}\tilde{\Lambda}_{21}+\tilde{\Lambda}_{21}^\top X_{12}^\top =-Y_{11}<0,\label{lyp21}\\
	&X_{11}\tilde{\Lambda}_{12}+X_{12}\tilde{\Lambda}_{22}=0.\label{lyp22}
\end{align}

{\color{red}Up to now, our objective is changed to prove the sufficient conditions for satisfying (\ref{lyp21})--(\ref{lyp22}) are conditions given in Theorem \ref{thm2}. To this end, we consider a system}
\begin{align}
	\dot{\zeta}=\tilde{\Lambda}_{22}\zeta.
\end{align}
{\color{red}A Lyapunov function can be chosen as $V(t)=\zeta^\top (R\otimes I_{n-l})\tilde{W}_{T,2}\zeta$ and its derivative along with $\zeta(t)$ gives rise to
\begin{align}
\dot{V}(t)=&\zeta^\top \tilde{W}_{T,2}diag\left\{r_i\begin{bmatrix}A_{22,io}-H_{i2o}C_{\diamond,i2o}&0\\A_{22,ir}&A_{22,iu}\end{bmatrix},\right.\notag\\
	&\quad\quad\quad\quad\quad\quad\quad\quad\quad\quad\quad\quad \left.i=1,\ldots,N\right\}\zeta,\notag\\
	&+\zeta^\top diag\left\{r_i\begin{bmatrix}A_{22,io}^\top -C_{\diamond,i2o}^\top H_{i2o}^\top &A_{22,ir}^\top \\0&A_{22,iu}^\top \end{bmatrix},\right.\notag\\
	&\quad\quad\quad\quad\quad\quad\quad\quad\quad\quad\quad\quad \left.
	i=1,\ldots,N\right\}\tilde{W}_{T,2}\zeta,\notag\\
	&-\gamma \zeta^\top T_{\diamond, 2}^\top ((R\mathcal{L}+\mathcal{L}^\top R)\otimes I_{n-l})T_{\diamond, 2}\zeta.
\end{align}
By using (\ref{w2}) and Lemma \ref{jointob2}, we know there is a $\mu>0$ such that
\begin{align}
	\dot{V}(t)=&-\gamma\zeta^\top  diag\left\{r_i\begin{bmatrix}I_{v_{i2}}&\\&0_{(n-l-v_{i2})\times n-l-v_{i2}}\end{bmatrix},\right.\notag\\
	&\quad\quad\quad\quad\quad\quad\quad\quad\quad\quad\quad\quad\quad\left.i=1,\ldots,N\right\}\zeta\notag\\
	&-\gamma\zeta^\top T_{\diamond, 2}^\top \left(\hat{\mathcal{L}}\otimes I_{n-l}\right)T_{\diamond, 2}\zeta+\zeta^\top (\tilde{\Lambda}_{22u}+\tilde{\Lambda}_{22u}^\top )\zeta\notag\\
	\leq&-\gamma\zeta^\top \left(T_{\diamond, 2}^\top \left(\hat{\mathcal{L}}\otimes I_{n-l}\right)T_{\diamond, 2}+G_{\diamond,2}\right)\zeta\notag\\
	&+\bar{\sigma}(\tilde{\Lambda}_{22u}+\tilde{\Lambda}_{22u}^\top )\|\zeta\|^2\notag\\
	\leq&-\gamma\mu\|\zeta\|^2+\bar{\sigma}(\tilde{\Lambda}_{22u}+\tilde{\Lambda}_{22u}^\top )\|\zeta\|^2,
\end{align}
}where
\begin{align}\label{lambda1}
	\tilde{\Lambda}_{22u}=&diag\left\{r_i\begin{bmatrix}0_{v_{i2}\times v_{i2}}&0_{v_{i2}\times(n-l-v_{i2})}\\A_{22,ir}&A_{22,iu}\end{bmatrix},\right.\notag\\
	&\quad\quad\quad\quad\quad\quad\quad\quad\quad\quad\quad\left.~i=1,\ldots,N\right\}.
\end{align}
Therefore, $\dot{V}<0$ if $\gamma>\lambda_1/\mu$ with $\lambda_1=\bar{\sigma}(\tilde{\Lambda}_{22u}+\tilde{\Lambda}_{22u}^\top )$. This indicates that $\tilde{\Lambda}_{22}$ is Hurwitz under the conditions of Theorem \ref{thm2}, and thus is invertible.

As a result, (\ref{lyp22}) holds if and only if $X_{12}=-\tilde{\Lambda}_{22}^{-1}X_{11}\tilde{\Lambda}_{12}$. Substitute $X_{12}$ into (\ref{lyp21}), and then we need to find $X_{11}>0$ such that
\begin{align}
	sym\{X_{11}\tilde{\Lambda}_{11}\}-\tilde{\Lambda}_{22}^{-1}X_{11}\tilde{\Lambda}_{12}\tilde{\Lambda}_{21}-\tilde{\Lambda}_{21}^\top \tilde{\Lambda}_{12}^\top X_{11}\tilde{\Lambda}_{22}^{-T}<0.\label{lyp23}
\end{align} 
{\color{red}Choose $X_{11}=(R\otimes I_l)\tilde{W}_{T,1}$ and use the results of Lemma \ref{jointob2}, then we obtain
\begin{align}
	&X_{11}\tilde{\Lambda}_{11}+\tilde{\Lambda}_{11}^\top X_{11}\notag\\
	=&X_{11}diag\left\{r_i\begin{bmatrix}A_{11,io}-H_{i1o}C_{\diamond,i1o}&\\&0_{(l-v_{i1})\times(l-v_{i1})}\end{bmatrix},\right.\notag\\
	&\quad\quad\quad\quad\quad\quad\quad\quad\quad\quad\quad\quad\quad\quad\quad\left.i=1,\ldots,N\right\},\notag\\
	&+diag\left\{r_i\begin{bmatrix}A_{11,io}^\top -C_{\diamond,i1o}^\top H_{i1o}^\top &\\&0_{(l-v_{i1})\times(l-v_{i1})}\end{bmatrix},\right.\notag\\
	&\quad\quad\quad\quad\quad\quad\quad\quad\quad\quad\quad\quad\quad\quad\quad\left.i=1,\ldots,N\right\}X_{11},\notag\\
	&+X_{11}\tilde{\Lambda}_{11u}+\tilde{\Lambda}_{11u}^\top X_{11}-\gamma X_{11}\tilde{W}_{T,1}^{-1}T_{\diamond 1}^\top (\mathcal{L}\otimes I_l)T_{\diamond 1}\notag\\
	&-\gamma T_{\diamond 1}^\top (\mathcal{L}^\top \otimes I_l)T_{\diamond 1}\tilde{W}_{T,1}^{-1}X_{11}\notag\\
	\leq&-\gamma \left(G_{\diamond,1}-T_{\diamond 1}^\top\left(\hat{\mathcal{L}\otimes I_l}\right)T_{\diamond 1}\right)+X_{11}\tilde{\Lambda}_{11u}+\tilde{\Lambda}_{11u}^\top X_{11}\notag\\
	<&-\gamma\mu I_{Nl}+X_{11}\tilde{\Lambda}_{11u}+\tilde{\Lambda}_{11u}^\top X_{11},\label{lyp24}
\end{align}
}where
\begin{align}
	\tilde{\Lambda}_{11u}=diag\left\{r_i\begin{bmatrix}0_{v_{i1}\times v_{i2}}&0_{v_{i2}\times(l-v_{i1})}\\A_{11,ir}&A_{11,iu}\end{bmatrix},~i=1,\ldots,N\right\}.
\end{align}
By combining (\ref{lyp23}) and (\ref{lyp24}), the sufficient conditions of (\ref{lyp23}) give rise to
\begin{align}
	&-\gamma\mu I_{Nl}+X_{11}\tilde{\Lambda}_{11u}+\tilde{\Lambda}_{11u}^\top X_{11}\notag\\
	&\quad-\tilde{\Lambda}_{22}^{-1}X_{11}\tilde{\Lambda}_{12}\tilde{\Lambda}_{21}-\tilde{\Lambda}_{21}^\top \tilde{\Lambda}_{12}^\top X_{11}\tilde{\Lambda}_{22}^{-T}<0.\label{lyp25}
\end{align}
Denote $\lambda_2$ the maximum eigenvalue of $X_{11}\tilde{\Lambda}_{11u}+\tilde{\Lambda}_{11u}^\top X_{11}+\tilde{\Lambda}_{22}^{-1}X_{11}\tilde{\Lambda}_{12}\tilde{\Lambda}_{21}-\tilde{\Lambda}_{21}^\top \tilde{\Lambda}_{12}^\top X_{11}\tilde{\Lambda}_{22}^{-T}$. Then, (\ref{lyp25}) holds if $\gamma>\lambda_2/\mu$.

{\color{red}Therefore, by choosing $X_{11}=(R\otimes I_l)\tilde{W}_{T,1}$, $X_{12}=-\tilde{\Lambda}_{22}^{-1}X_{11}\tilde{\Lambda}_{12}$, we can find $X_{22}>0$ such that
\begin{align*}
	X_{22}-X_{12}^\top X_{11}^{-1}X_{12}>0,
\end{align*} 
if $\gamma>\frac{1}{\mu}\max\{\lambda_1,\lambda_2\}$, which indicates $X>0$ owing to Schur complement lemma. Moreover, let $Y=diag\{(-\gamma\mu+\lambda_2)I_{l},I_{n-l}\}$ and thus (\ref{lyap2}) is satisfied based on the selected $X$ and $Y$. Therefore, error dynamics (\ref{xi-21})--(\ref{xi-22}) is admissible.} In other words, distributed observer (\ref{do2}) can achieve omniscience asymptotically and is impulse-free. 
\end{IEEEproof}

{\color{blue}Similar to Section \ref{sec3.1}, in order to facilitate readers in using the method of the theorem in this subsection, we will provide a brief guideline for designing parameters of (\ref{do2}).\\
	\textbf{STEP 1:}~According to the approach of singular value decomposition, there are two non-singular matrices $P_\diamond$ and $Q_\diamond$ such that $Q_\diamond EP_\diamond =diag\{	I_l,0_{(n-l)\times (n-l)}\}$. Then, we can obtain $A_{11}, A_{12}, A_{21}, A_{22}$ based on $Q_\diamond AP_\diamond$.\\
	\textbf{STEP 2:}~Calculate matrices related to observability decomposition based on (\ref{A3})--(\ref{A4}) of DDF.\\
	\textbf{STEP 3:}~Calculate observer gain $H_i$ based on the description in 2) of Theorem \ref{thm2}.\\
	\textbf{STEP 4:}~Weight matrices $\tilde{W}_{T,i1}$ and $\tilde{W}_{T,i2}$ can be solved by Lyapunov functions (\ref{w1}) and (\ref{w2}), respectively. Then, obtain $W_i$ in light of the description in 3) of Theorem \ref{thm2}.\\
	\textbf{STEP 5:}~Calculate coupling gain $\gamma$ according to 4) in Theorem 2).}

\begin{rem}
	\color{cyan}Note that step 5) in above guideline is complex and also requires the eigenvalues of $\hat{\mathcal{L}}$. Therefore, adaptive coupling gain $\gamma=\omega$ and $\dot{\omega}=\|\sum_{j=1}^N\alpha_{ij}(\hat{x}_j-\hat{x}_i)\|^2$
	is also recommended for distributed observer (\ref{do2}). 
\end{rem}

\begin{rem}\label{rem-assum3}
	\color{cyan}In a descriptor system, the observability of $(C_{\diamond 2}, A_{22})$ is equivalent to I-observability of the descriptor system. Therefore, we will now discuss the relationship between the observability of $(C_{\diamond 1},A_{11})$ and R-observability. In fact, R-observability is equivalent to $(C_{\diamond 1},A_{11}-A_{12}(A_{22}-H_{i2}C_{\diamond,i2})^{-1}A_{21})$ is observable. This means the observability of $(C_{\diamond 1},A_{11})$ is stronger than R-observability. However, based on this assumption, we no longer need the condition of $v_{i2}=n_2$ in subsection \ref{sec3.1}, so the distributed observer (\ref{do2}) is applicable for some systems that cannot meet the requirements of (\ref{do}). Lots of systems can satisfy Assumption \ref{assum3}, such as the hydraulic system in subsection \ref{sec4.1}  and the electronic network system given in subsection \ref{sec4.2}. 
\end{rem}

This subsection has designed the second type of distributed observer for descriptor systems based on Assumption \ref{assum3}. Combining with the conclusion in Subsection \ref{sec3.1}, the distributed observer method proposed in this paper can meet the requirements of a lot of descriptor systems, such as hydraulic system, mechanical system, isothermal reaction in an isothermal batch reactor system, economic system, and so on \cite{Belov2018,boubaker2019new}.

\subsection{Comparison between two designs}\label{sec3.3}
{\color{blue}
The previous two subsections presented two design methods for distributed observers of descriptor systems under different assumptions. In this subsection, we will compare these two different designs in two aspects.

1)~In terms of assumptions. Two different designs depend on two different assumptions. As stated in Remark \ref{rem-assum3}, the difference between Assumption \ref{assum2} and \ref{assum3} is that Assumption \ref{assum2} requires the system to be R-observable, while Assumption \ref{assum3} requires $(C_{\diamond 1},A_{11})$ to be observable. In general, the system is almost surely R-observable if $(C_{\diamond 1},A_{11})$ is observable. On the flip side, R-observable cannot deduce $(C_{\diamond 1},A_{11})$ is observable. For example, see the system with $E,A,C$ as follows:
\begin{align*}
	&E=diag\{I_3,0_{2\times 2}\},\\
	&A=\begin{bmatrix}
		0&0&0&1&-1\\a_1&a_2&a_3&0&1\\
		0&0&0&1&0\\-1&2&0&-1&0\\0&0&1&1&0		
	\end{bmatrix},~C=\begin{bmatrix}1&0&0&0&1\\0&0&0&0&1\\0&0&1&0&0\end{bmatrix}.
\end{align*}
It is easy to verify that this system is R-observable, but $(C_{\diamond 1},A_{11})$ is not observable. 

Therefore, in summary, Assumption \ref{assum3} is stronger than Assumption \ref{assum2}. However, the distributed observer based on Assumption \ref{assum2} requires an additional condition $v_{i2}=n$. As a result, both designs of distributed observers for descriptor systems can solve problems that the other one cannot do. They form complementarity in the application scope.

2)~In terms of parameter design. The distributed observer in subsection \ref{sec3.1} requires fewer design parameters (as seen in Theorem 1, only one pole placement and one Lyapunov equation need to be solved, while Theorem 2 requires two solutions), but the design of this distributed observer requires solving SDF, while implementation requires solving DDF. Therefore, the design process of the distributed observer in subsection \ref{sec3.1} is more complex. If a system can simultaneously satisfy the conditions of Theorem \ref{thm1} and Theorem \ref{thm2}, we recommend using the distributed observer in subsection \ref{sec3.2}.
}

\section{Implementation methods}\label{sec4}

{\color{cyan}Section \ref{sec3} has successfully completed the design and analysis of the descriptor distributed observer. However, it should be pointed out that the observers described by (\ref{do}) and (\ref{do2}) cannot be directly used to simulate the operation process of the underlying system (\ref{sys-11}). Observer, as a dynamic system, operates by spontaneously calculating the values in subsequent time as long as an initial value is given. However, descriptor systems are not normal dynamic systems and do not have such mechanisms. For example, in equation (\ref{ddf2}), we cannot calculate the value of $x_2(t)$ from $x_2(0)$ and $x_1(0)$ if $A_{22}$ is singular.
	
If such a situation exists in the observer dynamics, the observer cannot operate normally. Therefore, this section needs to analyze whether distributed observers (\ref{do}) and (\ref{do2}) can operate normally and provide the implementation methods of (\ref{do}) and (\ref{do2}) for their normal operation.}

\subsection{Implementation for (\ref{do})}\label{sec4.1}

The existence of impulse causes that (\ref{do}) cannot be directly operated. To overcome this issue, we first introduce $col\{\hat{x}_{i1},\hat{x}_{i2}\}=P_{\diamond}^{-1}\hat{x}_i$, and 
\begin{align}
	&Q_\diamond(A-H_iC_i)P_\diamond\triangleq\begin{bmatrix}A_{\diamond,i11}&A_{\diamond,i12}\\A_{\diamond,i21}&A_{\diamond,i22}\end{bmatrix},\notag\\
	&Q_\diamond W_i^{-1}P_\diamond\triangleq\begin{bmatrix}W_{\diamond,i11}&W_{\diamond,i12}\\W_{\diamond,i21}&W_{\diamond,i22}\end{bmatrix},~Q_\diamond H_i=\begin{bmatrix}
		H_{\diamond,i1}\\H_{\diamond,i2}
	\end{bmatrix}.\notag
\end{align}
Then, (\ref{do}) can be transformed into DDF, i.e.,
\begin{align}
	\dot{\hat{x}}_{i1}=&A_{\diamond,i11}\hat{x}_{i1}+A_{\diamond,i12}\hat{x}_{i2}+H_{i1}y_i\notag\\
	&+\gamma W_{\diamond,i11}\sum_{j=1}^N\alpha_{ij}\left(\hat{x}_{j1}-\hat{x}_{i1}\right),\label{imple-11}\\
	0=&A_{\diamond,i21}\hat{x}_{i1}+A_{\diamond,i22}\hat{x}_{i2}+H_{i2}y_i\notag\\
	&+\gamma W_{\diamond,i21}\sum_{j=1}^N\alpha_{ij}\left(\hat{x}_{ji}-\hat{x}_{i1}\right).\label{imple-12}
\end{align}
{\color{blue}The second lines of (\ref{imple-11}) and (\ref{imple-12}) is obtained by $Q_\diamond EP_\diamond=diag\{I_{l},0_{(n-l)\times(n-l)}\}$. Note that (\ref{imple-11}) and (\ref{imple-12}) are feasible if and only if $\hat{x}_{i2}$ can be solved by (\ref{imple-12}).

According to Theorem \ref{thm1}, we know error dynamics of (\ref{do}) is impulse-free. Moreover, only $A_{\diamond,i22}\hat{x}_{i2}$ in (\ref{imple-12}) is related to $\hat{x}_{i2}$. Consequently, (\ref{imple-11})--(\ref{imple-12}) is impulse-free if and only if $A_{\diamond,i22}$ is invertible. Therefore, we can obtain from (\ref{imple-12}) that}
\begin{align}
	\hat{x}_{i2}=&-A_{\diamond,i22}^{-1}(A_{21}-H_{i2}C_{\diamond,i1})\hat{x}_{i1}-A_{\diamond,i22}^{-1}H_{i2}y_i\notag\\
	&-\gamma A_{\diamond,i22}^{-1}W_{\diamond,21}\sum_{j=1}^N\alpha_{ij}\left(\hat{x}_{ji}-\hat{x}_{i1}\right).\label{imple-13}
\end{align}
Substituting (\ref{imple-13}) into (\ref{imple-11}) yields
\begin{align}
	\dot{x}_{i1}=&\left(A_{11}-A_{\diamond,i22}^{-1}(A_{21}-H_{i2}C_{\diamond,i1})\right)\hat{x}_{i1}\notag\\
	&-H_{i1}C_{\diamond,i1}\hat{x}_{i1}+\left(H_{i1}-A_{\diamond,i22}^{-1}H_{i2}\right)y_i\notag\\
	&+\gamma\left(W_{\diamond,11}-\gamma A_{\diamond,i22}^{-1}W_{\diamond,21}\right)\times\sum_{j=1}^N\alpha_{ij}\left(\hat{x}_{ji}-\hat{x}_{i1}\right).\label{imple-14}
\end{align}
Note that (\ref{imple-14}) is an ordinary differential equation (ODE) that can be operated by computer normally. In contrast (\ref{imple-13}) is an ordinary equation that can obtain $\hat{x}_{i2}$ by the knowledge of $\hat{x}_{i1}$ and $\hat{x}_{j1}$ from neighbors. Therefore, after solving the distributed observer parameters according to the method in Theorem \ref{thm1}, the normal operation of the distributed observer can be implemented by converting them into the forms of (\ref{imple-13}) and (\ref{imple-14}).

\subsection{Implementation for (\ref{do2})}\label{sec4.2}
We begin with transform (\ref{do2}) into
\begin{align}
	\dot{\hat{x}}_{i1}=&(A_{11}-H_{i1}C_{\diamond, i1})\hat{x}_{i1}+(A_{12}-H_{i1}C_{\diamond, i2})\hat{x}_{i2}+H_{i1}y_i\notag\\
	-&\gamma W_{T,i1}^{-1}\sum_{j=1}^N \alpha_{ij}(\hat{x}_{j1}-\hat{x}_{i1}),\\
	0=&(A_{21}-H_{i2}C_{\diamond, i1})\hat{x}_{i1}+(A_{22}-H_{i2}C_{\diamond, i2})\hat{x}_{i2}+H_{i2}y_i\notag\\
	-&\gamma W_{T,i2}^{-1}\sum_{j=1}^N \alpha_{ij}(\hat{x}_{j2}-\hat{x}_{i2}).\label{imple-17}
\end{align}
Different from (\ref{imple-12}), system (\ref{imple-17}) contains the coupling term with respect to $\hat{x}_{i2}$. As a result, the admissibility of error dynamics proved by Theorem \ref{thm2} cannot guarantee $A_{22}-H_{i2}C_{\diamond, i2}$ to be Hurwitz. However, there is a proper $\gamma$ such that $M_i\triangleq A_{22}-H_{i2}C_{\diamond,i2}-\gamma W_{T,i2}^{-1}\sum_{j=1}^N\alpha_{ij}I_{n-l}$ is invertible. Hence, we have
\begin{align}
	\dot{\hat{x}}_{i1}=&(A_{11}-H_{i1}C_{\diamond,i1})\hat{x}_{i1}+H_{i1}y_i\notag\\
	&+\gamma W_{T,i1}^{-1}\sum_{j=1}^N\alpha_{ij}\left(\hat{x}_{j1}-\hat{x}_{i1}\right)\notag\\
	&+(A_{12}-H_{i1}C_{\diamond,i2})M_{i}^{-1}\times (A_{21}-H_{i2}C_{\diamond,i1})\hat{x}_{i1}\notag\\
	&+(A_{12}-H_{i1}C_{\diamond,i2})M_i^{-1}H_{i2}y_i\notag\\
	&+(A_{12}-H_{i1}C_{\diamond,i2})M_i^{-1}W_{T,i2}^{-1}\sum_{j=1}^N\alpha_{ij}\hat{x}_{j2},\label{imple-15}\\
	\hat{x}_{i2}=&M_i^{-1}(A_{21}-H_{i2}C_{\diamond,i1})\hat{x}_{i1}+M_{i}^{-1}H_{i2}y_i\notag\\
	&+\gamma M_i^{-1}W_{T,i2}^{-1}\sum_{j=1}^N\alpha_{ij}\hat{x}_{j2}.\label{imple-16}
\end{align} 

Now, (\ref{imple-15}) and (\ref{imple-16}) can be operated normally. (\ref{imple-15}) is an ODE with external input $\hat{x}_{j1}$, $\hat{x}_{j2}$, and $y_i$. Since $\hat{x}_{j2}$ can be calculated by the neighbor agent through the similar equation as (\ref{imple-16}), $\hat{x}_{i1}$ can be obtained at any time. Then, $\hat{x}_{i2}$ can be calculated based on  $\hat{x}_{j1}$, $\hat{x}_{j2}$, and $y_i$. Therefore, we can implement distributed observer (\ref{do2}) by using (\ref{imple-15}) and (\ref{imple-16}).

\begin{rem}
	The centralized observers of the descriptor system also face the challenge of being unable to be directly operated. However, this issue can be easily overcome in centralized situations. The implementation method can be directly observed from the observer dynamics, while the implementation method of the descriptor system distributed observer is not intuitive. For this reason, we provide detailed implementation methods in this section.
\end{rem}

\subsection{Challenges in applying for real-world systems}
{\color{blue}
The first two subsections elaborate on the implementation and usage issues of the distributed observer for descriptor systems. In this subsection, we will briefly discuss our designed distributed observer's challenges in real-world systems and outline potential solutions.

1)~The regularity of the system. At present, most descriptor system theories are based on regular systems. The prerequisite for a system to be regular is that $E$ and $A$ must be square matrices. Moreover, even if $E$ and $A$ are square matrices, the system may not necessarily be regular, which poses certain challenges for applying the distributed observer designed in this paper in practical systems. Fortunately, under the premise that $E$ and $A$ are square matrices, if the system is not regular, then adding appropriate perturbations $\Delta A$ to the system matrix $A$ can make the pair $(E,A+\Delta A)$ regular. Example 2 in section \ref{sec4.2} is in this situation. However, adding perturbations can lead to model mismatch, which is also one of the challenges that needs to be addressed.

2)~Model mismatch. When describing the real-world system as a state space equation, the systems inevitably have model mismatches. In addition, to ensure the regularity of the systems, we sometimes need to introduce model mismatches artificially. The distributed observer designed in this paper cannot deal with model mismatches, but this issue is often addressed in existing theories of distributed observers for normal systems \cite{Xu2022TIV,XU2023ISA}. Therefore, to address the impact of model mismatch on our designed distributed observer, we can follow the normal system distributed observer theory and introduce a high-gain parameter in the observer gain $H_i$ when applying it to real-world systems. This method ensures that the error dynamics of the distributed observer can converge to any small compact set under the interference of model mismatch.

3)~Measurement noise. Measurement noise is an inevitable interference factor in real systems, and introducing high-gain parameters to address model mismatches will further amplify the negative impact of measurement noise. Distributed observer designed in this paper currently lacks the ability to cope effectively with measurement noise. However, introducing filtering and utilizing low-power high-gain theories \cite{ASTOLFI2018169} can effectively reduce the impact of noise on distributed observers.}

\begin{figure*}[!t]
	\centering
	\includegraphics[width=17cm]{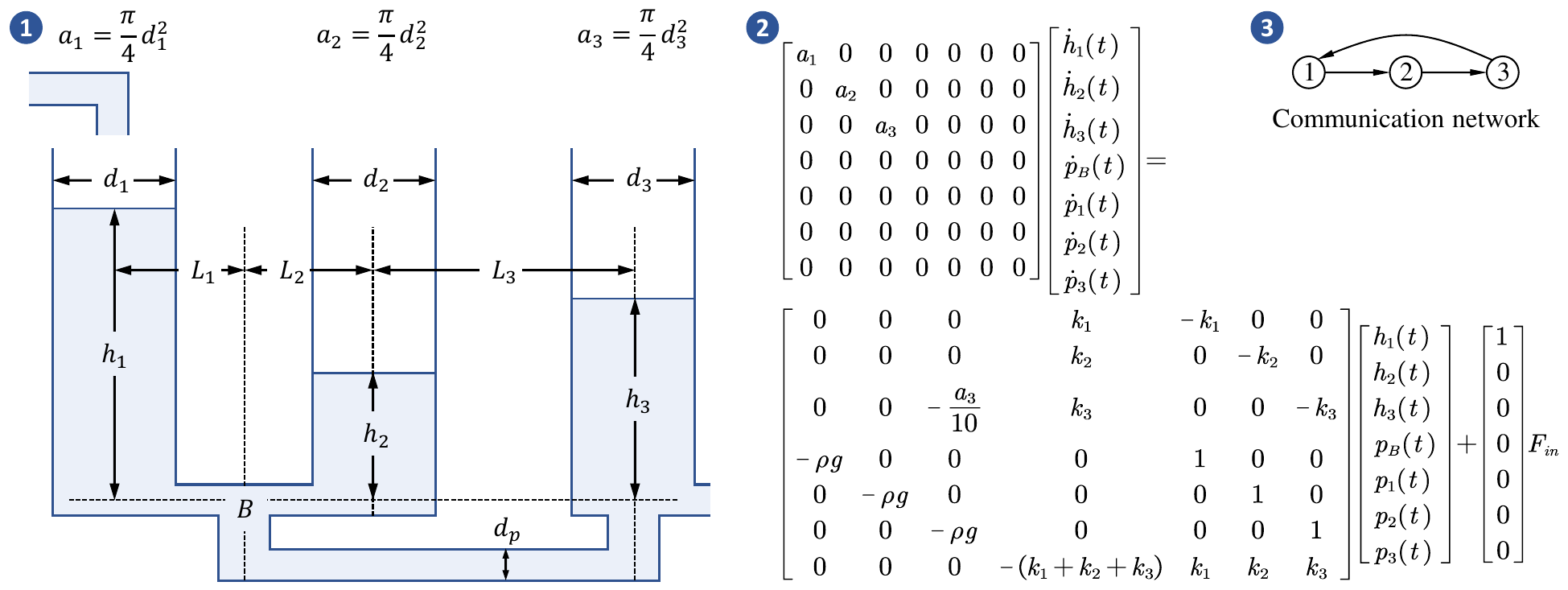}\\
	\caption{Hydraulic system and its dynamics described by descriptor system. Sub-figure 1 and sub-figure 2 are the construction and dynamics of hydraulic, respectively; sub-figure 3 shows the communication topology among three local observers.}\label{shuixiang}
\end{figure*}

\begin{figure*}[!t]
	\centering
	\includegraphics[width=17cm]{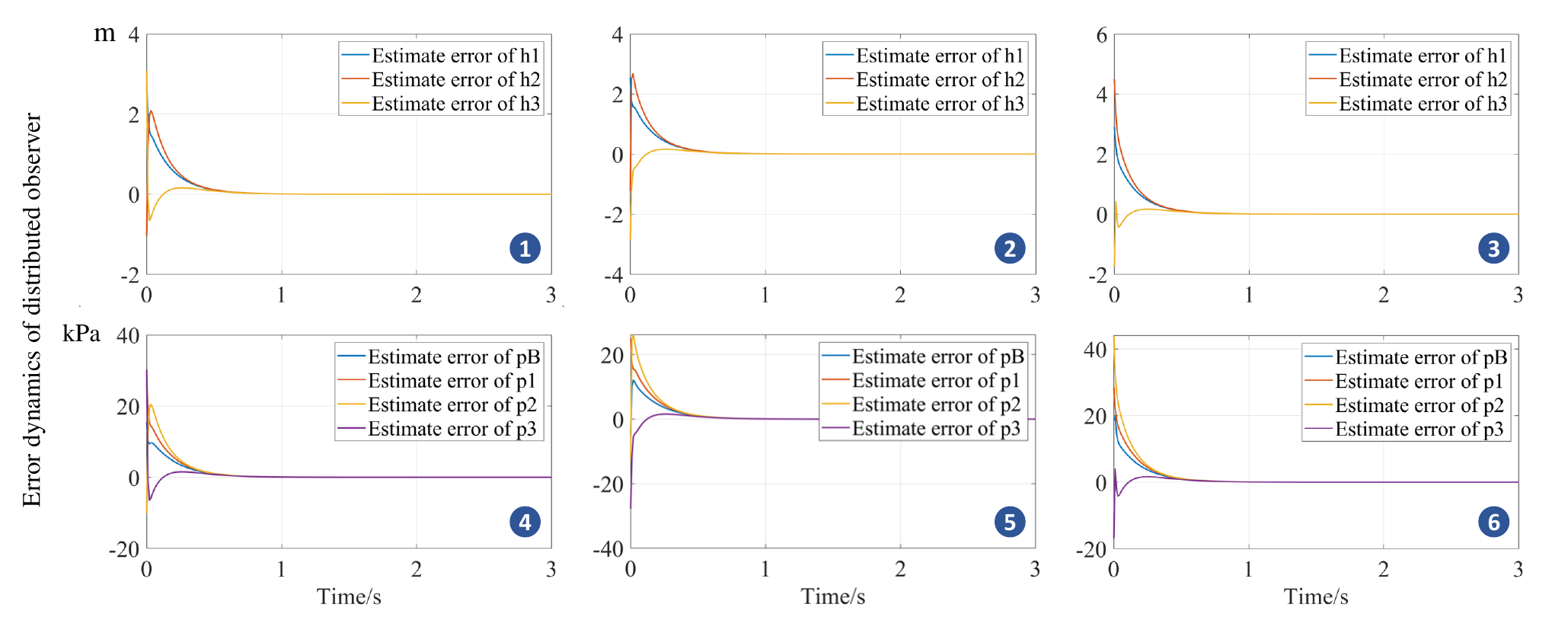}\\
	\caption{Error dynamics between distributed observer and real states. The first row and the second row shows the dynamics of $h_i$, and $p_i,p_B$, respectively, with $i=1,2,3$. Moreover, sub-figure 1, 4 is from local observer 1; sub-figure 2, 5 is from local observer 2; and sub-figure 3, 6 is from local observer 3.}\label{error-sx}
\end{figure*}

\begin{figure*}[!t]
	\centering
	\includegraphics[width=17cm]{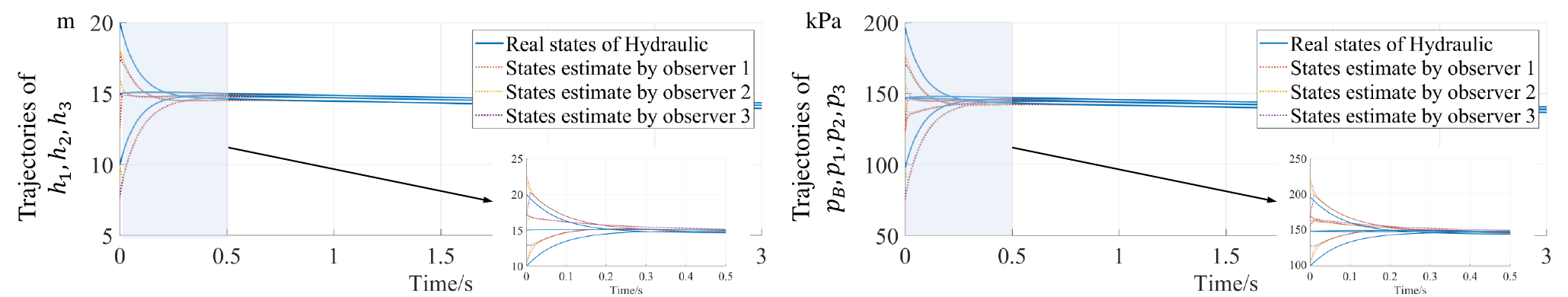}\\
	\caption{Trajectories of the real states and the state estimation of distributed observer. The solid line represents the true state, and the dashed line represents the state estimation of the state observer.}\label{states-sx}
\end{figure*}

\section{Simulations}\label{sec5}

{\color{cyan}In this section, two real-world examples are employed to display the effectiveness of the developed distributed observer. First, a hydraulic system which satisfies Assumption \ref{assum2} and $v_{i2}=n_2$ is given for verifying Theorem \ref{thm1} in Section \ref{sec5.1}. Moreover, to verify the results of Theorem \ref{thm2} as well as test the estimation ability of the proposed method for impulse phenomena, an electronic network system that satisfies Assumption 3 is considered in Section \ref{sec5.2}.}

\subsection{Simulation with hydraulic system}\label{sec5.1}

This subsection considers a model of the water height change in three tanks, an input flow in the first tank and the third tank leaking (See in the first sub-figure in Figure \ref{shuixiang}). 

The dimension of each part of the hydraulic system have been marked in Figure \ref{shuixiang}, including $d_1,d_2,d_3$, $L_1,L_2,L_3$, and $d_p$. The level of water in the three tanks are respectively $h_1,h_2,h_3$. The pressures at the bottom of these tanks are represented as $p_i=\rho gh_i$ for $i=1,2,3$, where $\rho=1.0\times 10^3{\rm kg/m^3}$ is the density of water and $g=9.8{\rm m/s^2}$ stands for the gravity acceleration. Let $p_B$ be the pressure at the pipe branch. Then, according to the Hagen-Poiseuille equation, the flow rates between the tanks and the pipe branch can be obtained:
\begin{align*}
	F_{B1}=&\left(p_B(t)-p_1(t)\right)\frac{\pi d_p^4}{128\eta L_B},\\
	F_{B2}=&\left(p_B(t)-p_2(t)\right)\frac{\pi d_p^4}{128\eta L_1},\\
	F_{B3}=&\left(p_B(t)-p_3(t)\right)\frac{\pi d_p^4}{128\eta L_2},
\end{align*}
where $\eta=10^{-3}{\rm Ns/m^2}$ is the dynamic viscosity of water. Since all water leaving tank 1 should enter into tank 2 and tank 3, we have $F_{B1}+F_{B2}+F_{B3}=0$. Let $F_{in}=1{\rm m^3/s}$ be the input of the hydraulic, then, the dynamics described by the descriptor system can be seen in the second sub-figure in Figure {\ref{shuixiang}}. In this system, $a_1,a_2,a_3$ are the cross-sectional areas of tanks 1, 2, and 3; and $k_1=\frac{\pi d_p^4}{128\eta L_B}$, $k_2=\frac{\pi d_p^4}{128\eta L_1}$, $k_3=\frac{\pi d_p^4}{128\eta L_2}$.

Consider three sensors that can measure $h_1,h_2,h_3$ respectively. Three local observers should estimate $h_i$, $p_i$, and $p_B$ with $i=1,2,3$ via the measurement outputs and the information exchange via the communication topology, which is shown in sub-figure 3 in Figure \ref{shuixiang}.

In this paper, we choose $d_1=d_2=d_3=10{\rm m}$, $d_p=2{\rm m}$, $L_B=L_1=10{\rm m}$, and $L_2=20{\rm m}$. By designing observer gains $H_i$ and weighted matrices $W_i$ in light of Theorem \ref{thm1}, we can implement the distributed observer with the methods in (\ref{do}) and (\ref{imple-15})--(\ref{imple-16}). Figure \ref{error-sx} shows the error dynamics of the distributed observer. It can be seen that all errors of all local observers can converge to zero asymptotically without impulse. Figure \ref{states-sx} compares the state estimation and the real states. In each sub-figure, the solid line represents the true state, and the dashed line represents the state estimation of the state observer. All of the above simulation results indicate the effectiveness of the proposed methods.

\begin{figure*}[!t]
	\centering
	\includegraphics[width=17cm]{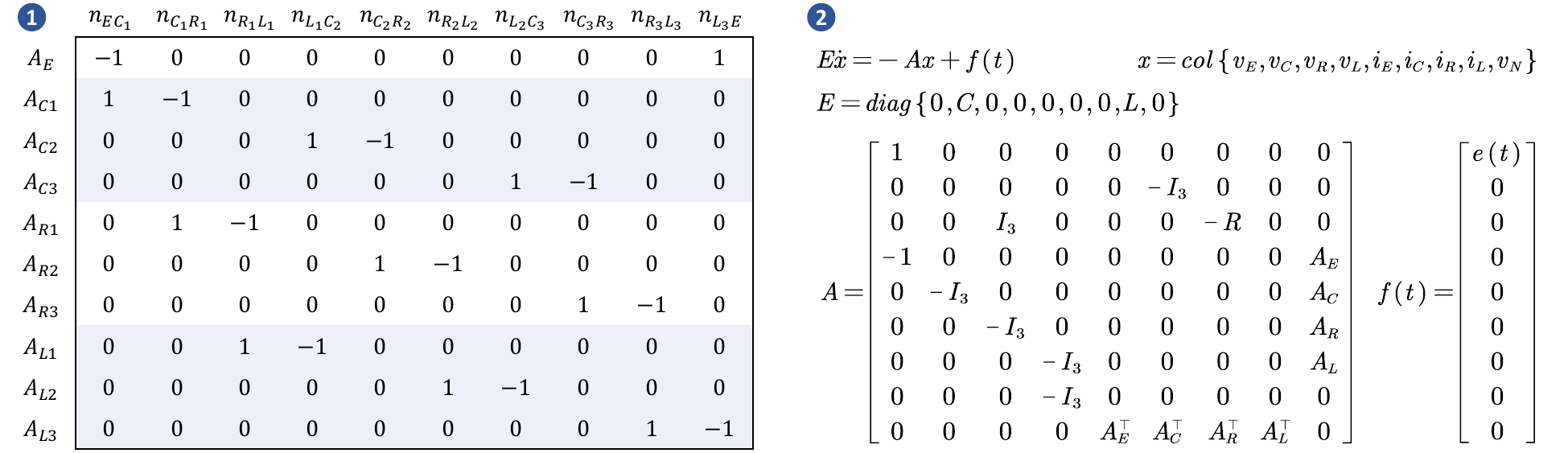}\\
	\caption{System formulation of electrical network system. Sub-figure 1 shows the incidence matrix, in which $\mathscr{C}_i$, $\mathscr{R}_i$, $\mathscr{L}_i$ stand for the $i$th capacitor, resistor, and inductor, respectively. $A_{\mathscr{E}}$ and $A_{\mathscr{K}}=col\{A_{\mathscr{K}_3},A_{\mathscr{K}_3},A_{\mathscr{K}_3}\}$ are the sub-matrices of $A_{cir}$, where $\mathscr{K}=\mathscr{E},\mathscr{C},\mathscr{R},\mathscr{L}$. $n_{\mathscr{X}\mathcal{Y}}$ represents the node connecting elements $\mathscr{X}$ and $\mathscr{Y}$.}\label{dianwang}
\end{figure*}

\begin{figure*}[!t]
	\centering
	\includegraphics[width=17cm]{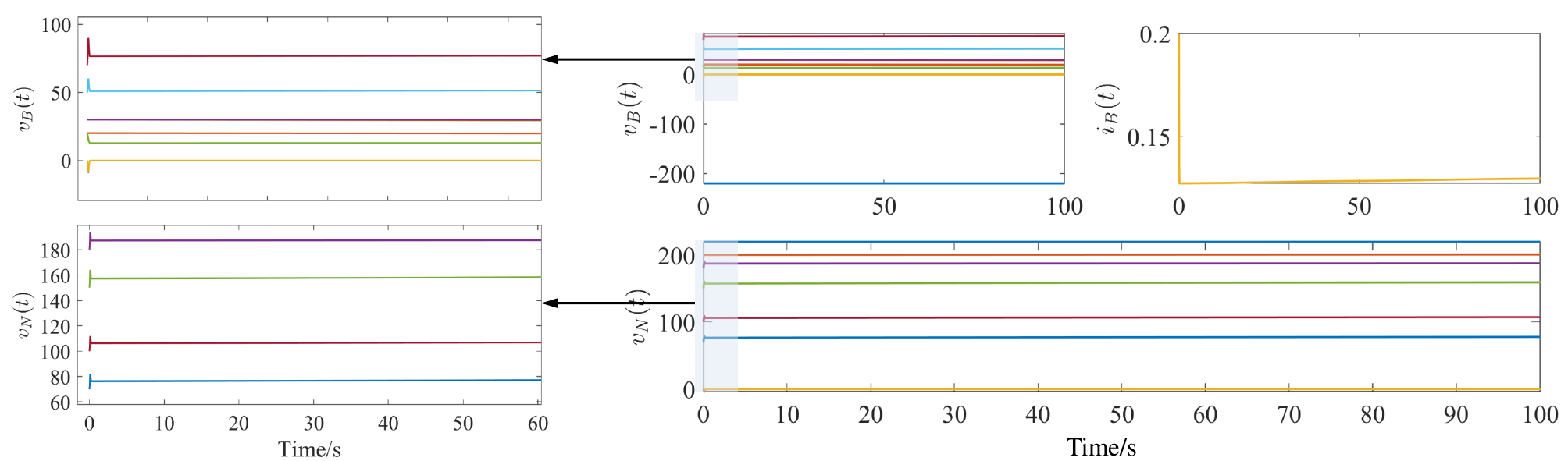}\\
	\caption{Trajectories of the states of electrical network system.}\label{error-dianwang}
\end{figure*}

\begin{figure*}[!t]
	\centering
	\includegraphics[width=17cm]{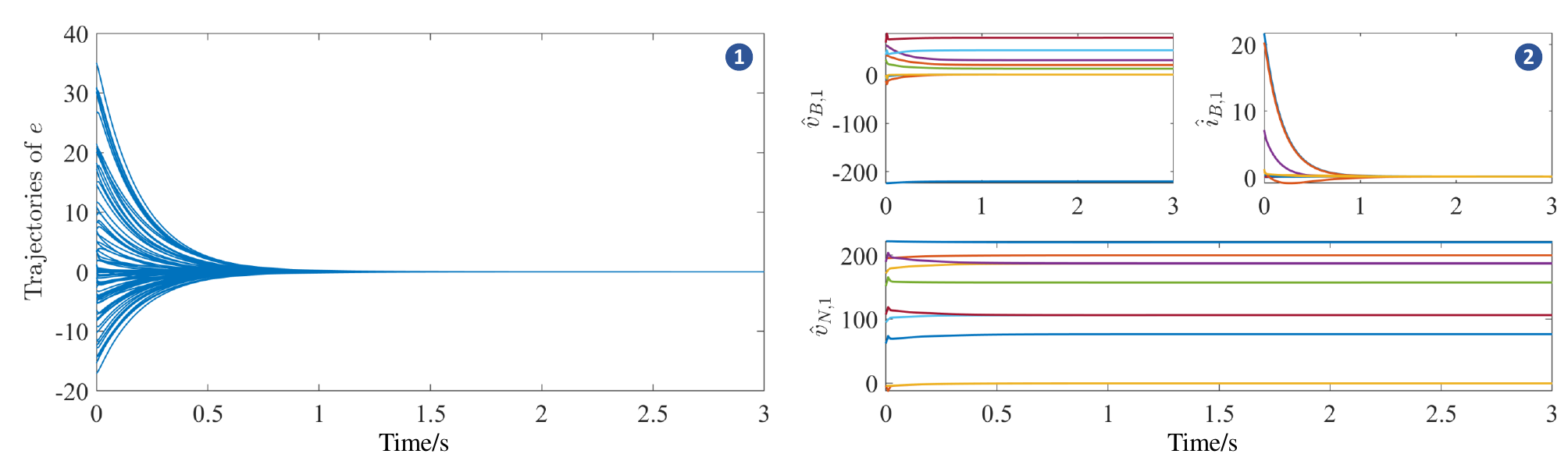}\\
	\caption{Error dynamics between distributed observer and real states.}\label{error-dianwang2}
\end{figure*}

\subsection{Simulation with electrical network}\label{sec5.2}

This subsection considers the electrical network consisting of $3$ branches, $1$ voltage source, and $10$ nodes. Each branch contains a resistor, a capacitor, and an inductor. Three local observers in this system are employed to estimate the state of the electrical network. The communication topology among them is the same as that in Example 1. Let $v_\mathscr{B}(t)=col\{v_\mathscr{E}(t),v_\mathscr{C}(t),v_\mathscr{R}(t),v_\mathscr{L}(t)\}$, $i_\mathscr{B}(t)=col\{i_\mathscr{E}(t),i_\mathscr{C}(t),i_\mathscr{R}(t),i_\mathscr{L}(t)\}$ be the branch voltages and branch currents respectively. The subscripts $\mathscr{R},\mathscr{C},\mathscr{L}\in\mathbb{R}^{3\times 3}$ are the diagonal matrices of resistors, capacitors, and inductors. $E=U\in\mathbb{R}$ is the electromotive force of the source with $U=220{\rm V}$. $v_\mathscr{K}$ and $i_\mathscr{K}$ stand for the voltage and current of the element $\mathscr{K}$ respectively, where $\mathscr{K}=\mathscr{E},\mathscr{C},\mathscr{R},\mathscr{L}$. The circuit is described by an incidence matrix $A_{cir}=col\{A_\mathscr{E},A_\mathscr{C},A_\mathscr{R},A_\mathscr{L}\}$ which is shown in the first sub-figure in Figure \ref{dianwang}. Denote $v_\mathscr{N}$ the nodal voltages and it satisfies $v_\mathscr{B}(t)=A_{cir}v_\mathscr{N}(t)$. According to Ohm's law and Kirchhoff's law, we have $v_\mathscr{R}(t)=i_\mathscr{R}(t)\mathscr{R}$ and $A_{cir}^\top i_\mathscr{R}(t)=0$. Moreover, there are linear capacitor and linear inductor voltage-current relations $i_\mathscr{C}(t)=\mathscr{C}\dot{v}_\mathscr{C}(t)$ and $v_\mathscr{L}(t)=\mathscr{L}\dot{i}_\mathscr{L}(t)$. Hence, this electrical network can be described by a descriptor system shown in the second sub-figure of Figure \ref{dianwang}.

In this system, the output matrix is given by
\begin{align*}
	C=\begin{bmatrix}
		1&0&0&0&0&0&0&0&0\\
		0&I_3&0&0&0&0&0&0&0\\
		0&0&I_3&0&0&0&0&0&0\\
		0&0&0&0&0&0&0&I_3&0\\
		0&0&0&0&0&0&0&0&I_{10}
	\end{bmatrix}.
\end{align*}
If $\mathcal{I}$ is an index set, we call $C_\mathcal{I}$ as the sub-matrix consists of the rows in $C$ indexed by $\mathcal{I}$. Hence, we can introduce the output measurements of each local observer as
\begin{align*}
	&C_1=C_{\mathcal{I}_1},~\mathcal{I}_1=\{1,2,5,8,11,12,13\},\\
	&C_2=C_{\mathcal{I}_2},~\mathcal{I}_2=\{3,6,9,14,15,16\},\\
	&C_3=C_{\mathcal{I}_3},~\mathcal{I}_3=\{4,7,10,17,18,19,20\}.
\end{align*} 

Set $\mathscr{C}=diag\{10,5,15\}{\rm \mu F}$, $\mathscr{R}=diag\{100,400,600\}{\rm \Omega}$, and $\mathscr{L}=diag\{1.2,1,1.1\}{\rm H}$. Now, the observer gains $H_i$, and weighted matrices $W_i$ can be calculated based on Theorem \ref{thm2} (Due to the high dimensionality of the parameter matrices, they will not be shown here). The initial states of the underlying system are chosen as $v_\mathscr{B}(0)=A_{cir}v_{\mathscr{N}}(0)$, $i_\mathscr{B}(0)=1_{10}\otimes 0.2{\rm A}$, and
\begin{align*}
	v_\mathscr{N}(0)=col\{220,200,180,180,150,100,100,70,0,0\}{\rm V}.
\end{align*} 
Then, the simulation results can be obtained and shown in Figure \ref{error-dianwang} and \ref{error-dianwang2} (Note that $A_{cir}^\top$ is not full row rank, so the system $(E,A)$ is not regular. To this issue, the entry at the bottom right of $A$ is set as $10^{-5}$ when simulating the system. Then, $(E,A)$ is regular, and the error caused by this entry can be ignored). 

Figure \ref{error-dianwang} shows the trajectories of the underlying system. As illustrated, it can be seen that this system has impulse at the initial time. Hence, this system can be employed to verify the effectiveness of the proposed method in estimating impulse.
Figure \ref{error-dianwang2} shows the effect of the designed distributed observer (only in the first 3 seconds). The first sub-figure shows the error dynamics of the distributed observer, whose trajectories do not contain impulses. It means the error dynamics are admissible. The second sub-figure displays the trajectories of the states estimate generated by local observer $1$ (States of other local observers are omitted because of the similar dynamics), which indicates that the developed distributed observer can effectively reconstruct the states of the descriptor system. 

\section{Conclusions and Future works}\label{sec6}
This paper have concerned distributed observer for descriptor system. The structure of distributed observer has been constructed by skillfully designing observer parameters. Then, sufficient conditions for asymptotic omniscience of distributed observer are given so that each local observer can reconstruct the states and impulse phenomenon of the underlying system. Moreover, distributed observer for normally linear system can be cashed as a particular case of this paper. Finally, simulation examples have shown that the designed distributed observer for descriptor system can achieve omniscience asymptotically. 

{\color{blue}In the future, we will further promote the results of this paper, including but not limited to extending it to nonlinear system situations using the theory of fully measured systems, extending it to large-scale systems based on our developed partitioned distributed observer theory, and extending it to dynamic time-varying topology situations using our proposed network transformation mapping \cite{Xu2025TAC}.}

\ifCLASSOPTIONcaptionsoff
  \newpage
\fi



%


\bibliographystyle{IEEEtran}
\bibliography{ref}

%

\begin{IEEEbiography}[{\includegraphics[width=1in,height=1.25in,clip,keepaspectratio]{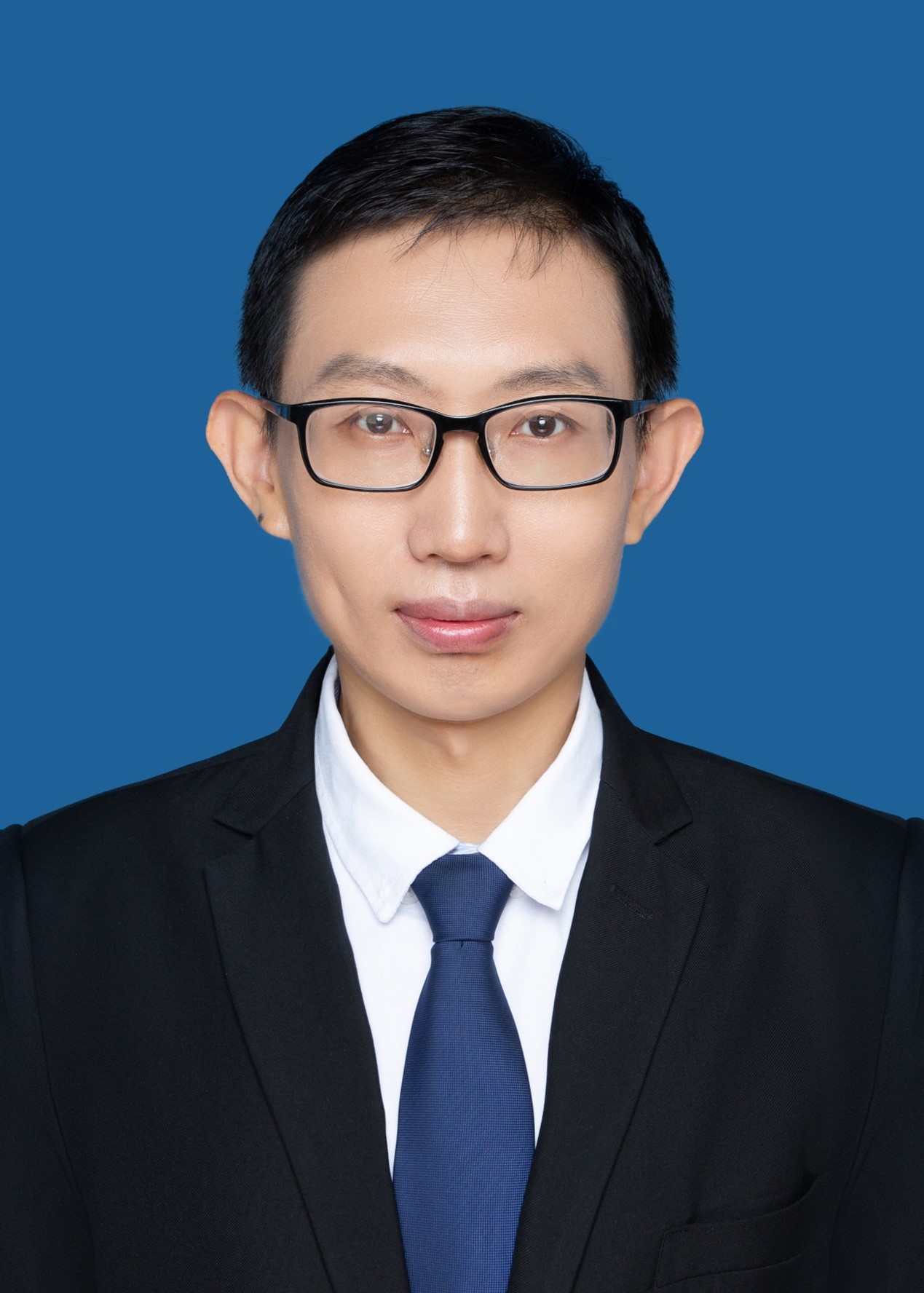}}]{Shuai Liu} (Member, IEEE) received the B.E. and M.E. degrees in control science and engineering from Shandong University, Jinan, China, in 2004 and 2007, respectively, and the Ph.D. degree in electrical and electronic engineering from Nanyang Technological University, Singapore, in 2012.
	
	From 2011 to 2017, he was a Senior Research Fellow with Berkeley Education Alliance, Singapore. Since 2017, he has been with the School of Control Science and Engineering, Shandong University. His research interests include distributed control, estimation and optimization, smart grid, integrated energy system, and machine learning. 
\end{IEEEbiography}

\begin{IEEEbiography}[{\includegraphics[width=1in,height=1.25in,clip,keepaspectratio]{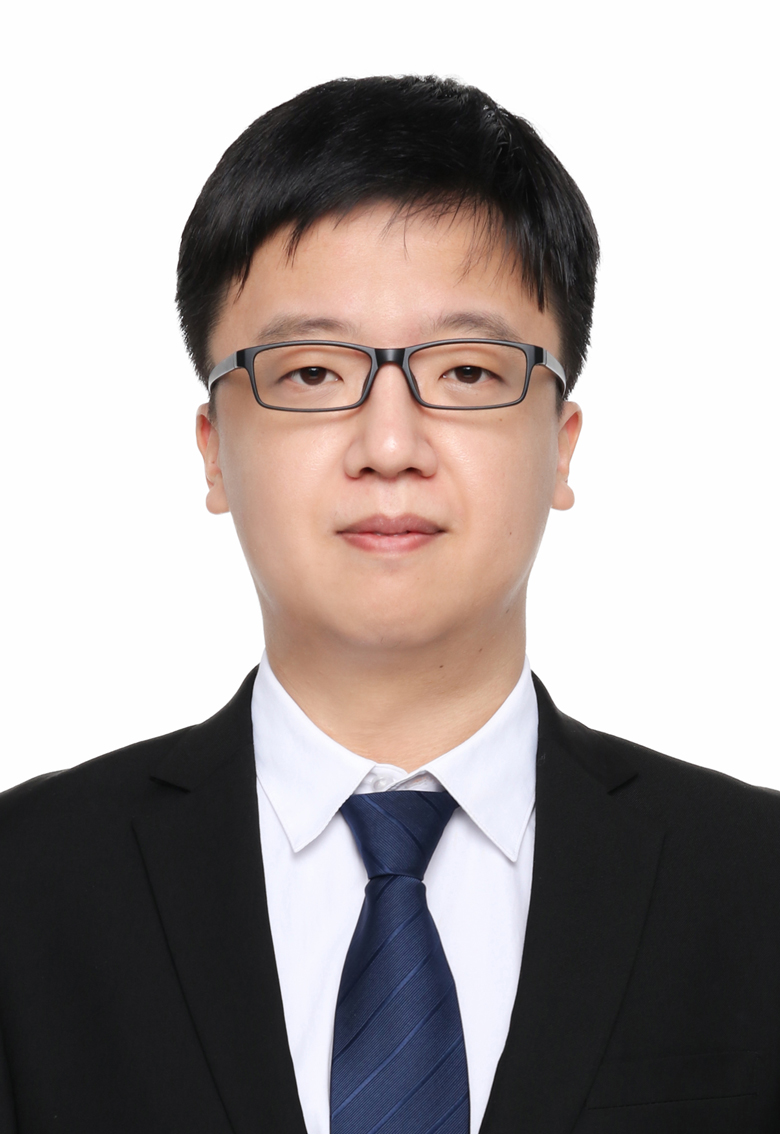}}]{Haotian Xu}
	received the B.S. degree in mathematics from Shandong University in 2016, and the Ph.D.degree in control science and engineering from Shanghai Jiao Tong University 2022. He is working as Post-Doctoral in control theory and engineering with Shandong University (2022–). His current research interests include distributed observers,cooperated stability based on distributed observer,nonlinear observer.
	
	He is the independent reviewer of international journals: IEEE Transactions on Automatic Control, IEEE Transactions on Neural Networks and Learning Systems, 
	IEEE Transactions on Systems, Man, and Cybernetics: Systems, IEEE Transactions on Industrial Electronics, IEEE Transactions on Control Systems Technology et al. His current research interests include Distributed observers, Cooperated stability based on Distributed observer, Nonlinear observer.
\end{IEEEbiography}
%
%




\end{document}